\documentclass[12pt, a4paper]{amsart}
\usepackage{fullpage}
\usepackage{amscd,amsmath,amsthm,amssymb,graphics,color}
\usepackage[colorlinks,linkcolor=blue,filecolor=blue,citecolor=blue,urlcolor=blue,bookmarks=false,hypertexnames]{hyperref} 
\usepackage{tikz}
\usepackage{caption,subcaption}

\def\NZQ{\mathbb}               
\def\NN{{\NZQ N}}
\def\QQ{{\NZQ Q}}
\def\ZZ{{\NZQ Z}}

\def\FF{{\NZQ F}}

\def\frk{\mathfrak}               

\def\mm{{\frk m}}

\def\Phi{{\frk n}}
\def\Phi{{\frk N}}

\def\MI{{\mathcal I}}

\def\MS{{\mathcal S}}

\def\ab{{\bold a}}

\def\opn#1#2{\def#1{\operatorname{#2}}} 
\opn\chara{char} \opn\length{\ell} \opn\pd{pd} \opn\rk{rk}
\opn\projdim{proj\,dim} \opn\injdim{inj\,dim} \opn\rank{rank}
\opn\depth{depth} \opn\grade{grade} \opn\height{height}\opn\coheight{coheight}
\opn\embdim{emb\,dim} \opn\codim{codim}

\opn\Tr{Tr} \opn\bigrank{big\,rank}
\opn\superheight{superheight}\opn\lcm{lcm}
\opn\trdeg{tr\,deg}
\opn\reg{reg} \opn\lreg{lreg} \opn\ini{in} \opn\lpd{lpd}
\opn\size{size}\opn\bigsize{bigsize}
\opn\cosize{cosize}\opn\bigcosize{bigcosize}
\opn\sdepth{sdepth}\opn\sreg{sreg}
\opn\link{link}\opn\fdepth{fdepth}\opn\type{type}
\opn\GL{GL} 

\opn\div{div} \opn\Div{Div} \opn\cl{cl} \opn\Cl{Cl}
\opn\Spec{Spec} \opn\Supp{Supp} \opn\supp{supp} \opn\Sing{Sing}
\opn\Ass{Ass} \opn\Min{Min}\opn\Mon{Mon} \opn\dstab{dstab} \opn\astab{astab}
\opn\Syz{Syz}
\opn\Ann{Ann} \opn\Rad{Rad} \opn\Soc{Soc}
\opn\Im{Im} \opn\ker{ker} \opn\Coker{Coker} \opn\Am{Am}
\opn\Hom{Hom} \opn\Tor{Tor} \opn\Ext{Ext} \opn\End{End}
\opn\Aut{Aut} \opn\id{id}

\opn\nat{nat}
\opn\pff{pf}
\opn\Pf{Pf} \opn\GL{GL} \opn\SL{SL} \opn\mod{mod} \opn\ord{ord}
\opn\Gin{Gin} \opn\Hilb{Hilb}\opn\sort{sort}
\opn\Proj{Proj}
\opn\aff{aff} \opn\con{conv} \opn\relint{relint} \opn\st{st}
\opn\lk{lk} \opn\cn{cn} \opn\core{core} \opn\vol{vol}
\opn\link{link} \opn\star{star}\opn\lex{lex}
\opn\gr{gr}

\opn\dirlim{\underrightarrow{\lim}}
\opn\inivlim{\underleftarrow{\lim}}
\let\dirsum=\oplus
\let\tensor=\otimes
\let\iso=\cong

\let\Dirsum=\bigoplus

\let\ab\allowbreak
\let\to=\rightarrow

\def\Implies{\ifmmode\Longrightarrow \else
	\unskip${}\Longrightarrow{}$\ignorespaces\fi}
\def\implies{\ifmmode\Rightarrow \else
	\unskip${}\Rightarrow{}$\ignorespaces\fi}
\def\iff{\ifmmode\Longleftrightarrow \else
	\unskip${}\Longleftrightarrow{}$\ignorespaces\fi}

\let\:=\colon
\newtheorem{Theorem}{Theorem}[section]
\newtheorem{Lemma}[Theorem]{Lemma}
\newtheorem{Corollary}[Theorem]{Corollary}
\newtheorem{Proposition}[Theorem]{Proposition}
\newtheorem{Remark}[Theorem]{Remark}

\newtheorem{Definition}[Theorem]{Definition}

\newtheorem{Notation}[Theorem]{Notation}

\newtheorem{Observation}[Theorem]{Observation}

\let\epsilon\varepsilon

\title {Green-Lazarsfeld property $N_p$ for Segre product of Hibi rings}

\author {Dharm Veer}

\dedicatory{In memory of Prof. C.~S.~Seshadri}
\address{Chennai Mathematical Institute, Siruseri, Tamilnadu 603103, India.}
\email{dharm@cmi.ac.in}

 \thanks{The author was partly supported by the grant CRG/2018/001592 (Manoj Kummini) from Science and Engineering Research Board, India and by an Infosys Foundation fellowship.}

\begin{document}

\begin{abstract}
In this article, we prove that if a Hibi ring satisfies property $N_2$, then its Segre product with a polynomial ring in finitely many variables also satisfies property $N_2$. When the polynomial ring is in two variables, we also prove the above statement for $N_3$. Moreover, we study the minimal Koszul relations of the second syzygy module of Hibi rings. 
\end{abstract} 

\maketitle

\section{Introduction}

A classical problem in commutative algebra is to study the graded minimal free resolution of graded modules over polynomial rings. Let $S$ be a standard graded polynomial ring in finitely many variables over a field $K$ and $I$ be a graded $S$-ideal. To study the graded minimal free resolution of $S/I$, Green-Lazarsfeld \cite{[GL86]} defined property $N_p$ for $p\in \NN$. The ring $S/I$ satisfies property $N_p$ if $S/I$ is normal and the graded minimal free resolution of $S/I$ over $S$ is linear upto $p^{th}$ position. In this article, we study the Green-Lazarsfeld property $N_p$ for Segre product of Hibi rings and the minimal Koszul sygygies of the Hibi ideals and the initial Hibi ideals.

Let $L$ be a finite distributive lattice and $P=\{p_1,\ldots,p_{n}\}$ be the subposet of join-irreducible elements of $L$. Let $K$ be a field and let $R=K[y_1,\ldots,y_n,z_1,\ldots,z_n]$ be a polynomial ring over $K$. The Hibi ring associated with $L$, denoted by $R[L]$,  is the subring of $R$ generated by the  monomials $u_{\alpha}=(\prod_{p_i\in \alpha}y_i)(\prod_{p_i\notin \alpha}z_i)$ where $\alpha\in L$. Hibi \cite{[HIBI87]} showed that $R[L]$ is a normal Cohen–Macaulay domain of dimension $\#P+1$, where $\#P$ is the cardinality of $P$. Let $K[L]=K[\{x_\alpha: \alpha\in L\}]$ be the polynomial ring over $K$ and $\pi: K[L] \to R[L]$ be the $K$-algebra homomorphism with $x_\alpha \mapsto u_{\alpha}$. The ideal $I_L=(x_\alpha x_\beta-x_{\alpha\wedge\beta}x_{\alpha\vee\beta}: \alpha,\beta\in L\ \text{and}\ \alpha,\beta \text{ incomparable})$ is the kernel of the map $\pi$. It is called the Hibi ideal associated to $L$.\par

In past, various authors have studied minimal free resolution of Hibi rings.
Ene et al.~\cite{[EHM15]} provided a combinatorial formula for the regularity of Hibi rings.
In~\cite{EQR13regularityplanardistlattices, [EHH15]}, the authors have characterized all Hibi rings with linear resolution.
Ene \cite{[ENE15]} characterizes all simple planar distributive lattices for which the associated Hibi ring satisfies property $N_2$. 
Das and Mukherjee~\cite{[DM17]} described the generators of second syzygy module of simple planar distributive lattices.
Additionally, the author~\cite{VEE21} has proved several results regarding property $N_p$ of Hibi rings, including a characterization of those Hibi rings that satisfy property $N_p$ for $p\geq 4$.

The Segre product of polynomial rings may be viewed as a Hibi ring. Let $A = K[x_{1,0},\ldots,x_{1,n_1}]* \cdots * K[x_{r,0},\ldots,x_{r,n_r}]$ be the Segre product of $r$ polynomial rings, where $n_i \geq 1$ and $n_i\in \NN$ for all $i$. Sharpe \cite{[SHARPE64]} proved in 1964 that if $r=2$, then $A$ satisfies property $N_2$. For $r=2$, Lascoux \cite{[LAS78]} and Pragacz-Weyman \cite{[PW85]} proved that $A$ satisfies property $N_3$ if $K$ contains the rational field $\QQ$. Hashimoto \cite{[HAS90]} showed that if $r=2$, $n_1,n_2\geq 4$ and  characteristic of the field $K$ is 3, then $A$ does not satisfy property $N_3$. Rubei \cite{[RUBEI02],[RUBEI07]} proved that if $r\geq 3$ and char$(K)=0$, then $A$ satisfies property $N_3$ but it does not satisfy property $N_4$. 

It is known that if a poset in disconnected, then the associated Hibi ring is the Segre product of the Hibi rings associated to each individual connected piece of the poset.
Let $P$ be a poset that is disjoint union of a poset $P_1$ and an isolated point. 
In this article, we prove that if $R[\MI(P_1)]$ satisfies property $N_p$ for $p\leq3$, then so does $R[\MI(P)]$.
Our proof follows the argument of Rubei \cite{[RUBEI02]} which is combinatorial in nature.
The core idea of Rubei's approach is as follows: 
in~\cite{BH, [STURM96]}, it was shown that multigraded Betti numbers of affine semigroup rings can be computed using homologies of squarefree divisor complexes.
Utilizing the results of~\cite{BH, [STURM96]}, 
it is suffices to prove that 
for $p=2$, the first homology of certain squarefree divisor complexes vanishes, 
whereas for $p=3$, the second homology of certain squarefree divisor complexes vanishes. 
For $p=2$ (resp. $p=3$), we begin by showing that every $1$-cycle (resp. $2$-cycle) $\gamma$ of the simplicial complex $\Delta$ is homologous to an $1$-cycle (resp. a $2$-cycle) of a subcomplex of $\Delta$. 
We then use this result to demonstrate that $\gamma$ is, in fact, a boundary of $\Delta$.

We generalize the above result for $p=2$ to the case when $P$ is the disjoint union of a poset $P_1$ and a chain. 
More precisely, we show that if $R[\MI(P_1)]$ satisfies property $N_2$, then $R[\MI(P)]$ also satisfies property $N_2$.
To prove this, we use induction on the length of the chain. 
The above result serves as the base case for the induction.
This proof is motivated from Rubei~\cite{[RUBEI07]}.
Following the combinatorial argument outlined in the previous result,
we first show that every $1$-cycle $\gamma$ of the simplicial complex $\Delta$ is homologous to an $1$-cycle $\gamma'$ of a simplicial complex of $X$, 
which contains $\Delta$ as a subcomplex. 
We then note that $\gamma'$ is homologous to $0$ is $X$ using the hypothesis $R[\MI(P_1)]$ satisfies property $N_2$.
Finally, we consider different cases to show that $\gamma$ is homologous to $0$ in $\Delta$, which is equivalent to saying that $\gamma$ is a boundary of $\Delta$.

Next, we inquire the following: if a Hibi ring does not satisfy property $N_2$, then which Koszul relations will be in the minimal generating set of its second syzygy module.
We provide a partial answer this question. In the process, we characterize the Koszul relations pairs of the initial Hibi ideals under a fixed monomial order.

The article is organized as follows. 
In Section~\ref{hibiringdefi}, we recall some basic notions of algebra and combinatorics.
Section~\ref{segre} is about the results on property $N_p$ of Segre products of Hibi rings.
Finally, in Section~\ref{koszulpairs}, we study the Koszul relation pairs of the Hibi ideals and initial Hibi ideals.

\subsection*{Acknowledgments}
I am extremely grateful to Manoj Kummini for his guidance and various insightful discussions throughout the preparation of this article. The computer algebra systems Macaulay2~\cite{[M2]} and SageMath~\cite{[sagemath]} provided valuable assistance in studying examples.

\section{Preliminaries}\label{hibiringdefi}
 
We start by defining some notions of posets and distributive lattices. For more details and examples, we refer the reader to \cite[Chapter\ 3]{[STAN12]}. Throughout this article, all posets and distributive lattices will be finite.
 
Let $P$ be a poset. We say that two elements $x$ and $y$ of $P$ are {\em comparable}  if $x\leq y$ or $y\leq x$; otherwise $x$ and $y$ are {\em incomparable}. For $x,y \in P$, we say that {\em y covers x} if $x<y$ and there is no $z\in P$ with $x<z<y$. We denote it by $x\lessdot y$. A poset is completely determined by its cover relations. A {\em chain} $C$ of $P$ is a totally ordered subset of $P$. The {\em length} of a chain $C$ of $P$ is $\#C - 1$.
A subset $\alpha$ of $P$ is called an {\em order ideal} of $P$ if it satisfies the following condition: for any $x \in \alpha$ and $y \in P$, if $y\leq x$, then  $y\in \alpha$. Define $\MI(P):= \{\alpha \subseteq P: \alpha  \ \text{is an order ideal of}\ P \}$. It is easy to see that $\MI(P)$, ordered by inclusion, is a distributive lattice under union and intersection. $\MI(P)$ is called the {\em ideal  lattice} of the poset $P$.  

Let $P$ and $Q$ be two posets. The {\em ordinal sum} $P\dirsum Q$ of the disjoint posets $P$ and $Q$ is the poset on the set $P\cup Q$ with the following order: if $x, y \in P\dirsum Q $, then $x \leq y$ if either $x, y \in P$ and $x \leq y$ in $P$ or $x, y \in Q$ and $x \leq y$ in $Q$ or $x \in P$ and $y \in Q$. Let $P,Q$ be posets on disjoint sets. The {\em disjoint union} of posets $P$ and $Q$ is the poset $P+Q$ on the set $P\cup Q$ with the following order: if $x,y \in P+Q$, then $x\leq y$ if either $x, y \in P$ and $x \leq y$ in $P$ or $x, y \in Q$ and $x \leq y$ in $Q$. A poset which can be written as direct sum of two posets is called {\em disconnected}. Otherwise, $P$ is {\em connected}.	
 
Let $L$ be a distributive lattice. An element $x\in L$ is called {\em join-irreducible} if $x$ is not the minimal element of $L$ and whenever $x =y \vee z$ for some $y,z \in L$, we have either $x=y$ or $x=z$. Let $P$ be the subposet of join-irreducible elements of $L$. By Birkhoff's fundamental structure theorem \cite[Theorem\ 3.4.1]{[STAN12]}, $L$ is isomorphic to the ideal lattice $\MI(P)$. Write $P=\{p_1,\ldots,p_n\}$ and let $R=K[y_1,\ldots,y_n,z_1,\ldots,z_n]$ be a polynomial ring in $2n$ variables over a field $K.$ The {\em Hibi ring} associated with $L$, denoted by $R[L]$,  is the subring of $R$ generated by the  monomials $u_{\alpha}=(\prod_{p_i\in \alpha}y_i)(\prod_{p_i\notin \alpha}z_i)$ where $\alpha\in L$.

Let $K[L]=K[\{x_\alpha: \alpha\in L\}]$ be the polynomial ring over $K$ and $\pi: K[L] \to R[L]$ be the $K$-algebra homomorphism with $x_\alpha \mapsto u_{\alpha}$.  Let $I_L=(x_\alpha x_\beta-x_{\alpha\wedge\beta}x_{\alpha\vee\beta}: \alpha,\beta\in L\ \text{and}\ \alpha,\beta \text{ incomparable})$ be an $K[L]$-ideal.

Let $<$ be a total order on the variables of $K[L]$ with the property that one has $x_\alpha < x_\beta$ if $\alpha < \beta$ in $L$. Consider the graded reverse lexicographic order $<$ on $K[L]$ induced by this order of the variables. 

\begin{Theorem}\cite[Theorem\ 6.19]{[HHO18]}\label{grobner}
	The generators of $I_L$ described above forms a Gr\"obner basis of $\ker(\pi)$ with respect to $<$. In particular, $\ker(\pi) = I_L$.
\end{Theorem}

The ideal $I_L$ is called the {\em Hibi ideal} of $L$. By Theorem~\ref{grobner}, it follows that 
\[\ini_<(I_L)=(x_\alpha x_\beta: \alpha,\beta\in L \ \text{and}\ \alpha, \beta \text{ incomparable}).\]\par

It is easy to see that if $\mathcal{P}=\{p_1,\ldots, p_{n}\}$ is a chain with $p_1\lessdot \cdots \lessdot p_n$, then $\MI(\mathcal{P}) =\{\emptyset,\{p_1\},\{p_1,p_2\},\ldots,\{p_1,\ldots,p_n\}\}$. Thus, the Hibi ring $R[\MI(\mathcal{P})] = K[z_1\cdots z_n,y_1z_2\cdots z_n,\ab \ldots,y_1\cdots y_n]$, which is a polynomial ring in $n+1$ variables. Let $P_1$ and $P_2$ be two posets and $P$ be their disjoint union. It was observed in \cite{[HHR00]} that  $R[\MI(P)] \cong R[\MI(P_1)] * R[\MI(P_2)]$, where $*$ denotes the Segre product. 

\subsection{Green-Lazarsfeld property}\label{propertynp} 
 Let $S=K[x_1,\ldots, x_n]$ be a standard graded polynomial ring in $n$ variables over a field $K$ and $M$ be a graded $S$-module. Let $\FF$ be the graded minimal free resolution of $M$ over $S$:
$$\FF : 0 \to \mathop{\Dirsum}_{j}^{} S(-j)^{\beta_{rj}}\to \cdots \to \mathop{\Dirsum}_{j}^{} S(-j)^{\beta_{1j}} \to \mathop{\Dirsum}_{j}^{} S(-j)^{\beta_{0j}}.$$

The numbers $\beta_{ij}$ are called the minimal graded Betti numbers of the module $M$.\par

Let $I$ be a graded $S$-ideal. We say that $S/I$ satisfies {\em Green-Lazarsfeld property $N_p$} if $S/I$ is normal and $\beta_{ij}(S/I)=0$ for $i \neq j+1$ and $1\leq i\leq  p$. Therefore, $S/I$ satisfies property $N_0$ if and only if it is normal; it satisfies property $N_1$ if and only if it is normal  and $I$ is generated by quadratics; it satisfies property $N_2$ if and only if it satisfies property $N_1$ and $I$ is linearly presented and so on. We know that the Hibi rings are normal and the Hibi ideals are generated by quadratics. Hence, the Hibi rings satisfy property $N_1$. 

\subsection{Squarefree divisor complexes}\label{squarefree}
Let  $H\subset \NN^n$ be an affine semigroup and $K[H]$ be the semigroup ring attached to it. Suppose that $h_1,\ldots, h_m\in \NN^n$ is the unique minimal set of generators of $H$. We consider the polynomial ring $T=K[t_1,\ldots,t_n]$ in $n$ variables. Then $K[H]$ is the subring of $T$ generated by the monomials $u_i=\prod_{j=1}^nt_j^{h_i(j)}$ for $1\leq i \leq m$, where $h_i(j)$ denotes the $j$th component of the integer vector $h_i$. Consider a $K$-algebra map $S=K[x_1,\ldots,x_m]\to K[H]$ with $x_i\mapsto u_i$ for $i=1,\ldots,m$.  Let $I_H$ be the kernel of this $K$-algebra map. Set  $\deg x_i=h_i$ to assign a $\ZZ^n$-graded ring structure to $S$. Let $\mm$ be the graded maximal $S$-ideal. Then $K[H]$ as well as $I_H$ become $\ZZ^n$-graded $S$-modules. Thus, $K[H]$ admits a minimal $\ZZ^n$-graded  $S$-resolution $\FF$.
  
Given $h\in H$, we define the {\em squarefree divisor complex} $\Delta_h$ as follows:

 $$\Delta_h:=\{F \subseteq [m] :\prod_{i\in F}^{} u_i \ \text{divides}\ t_1^{h(1)}\cdots t_n^{h(n)}\ \text{in}\ K[H]\}.$$

 We denote the $i$th reduced simplicial homology of a simplicial complex $\Delta$ with coefficients in $K$ by $\widetilde{H}_{i}(\Delta, K)$.
  
  \begin{Proposition}\cite[Proposition\ 1.1]{BH}, \cite[Theorem\ 12.12]{[STURM96]}
  	\label{bh}
  	With the notation and assumptions introduced one has  $\Tor_i(K[H],K)_h\iso\widetilde{H}_{i-1}(\Delta_h, K)$. In particular,
  	$$\beta_{ih}(K[H])=\dim_K\widetilde{H}_{i-1}(\Delta_h, K).$$
  \end{Proposition}
  
  \begin{Definition}
   Let $H\subset \NN^n$ be an affine semigroup generated by $h_1,\ldots, h_m$. An affine subsemigroup $H'\subset H$ generated by  a subset of $\{h_1,\ldots, h_m\}$ will be called a {\em homologically pure} subsemigroup of $H$ if for all $h\in H'$ and all $h_i$ with $h-h_i\in H$,  it follows that $h_i\in H'$.
  \end{Definition}
  
  Let $H'$ be a subsemigroup of $H$ generated by a subset $\mathcal{X}$ of $\{h_1,\ldots,h_m\}$, and let $S'=K[\{x_i : h_i\in \mathcal{X}\}]\subseteq S$. Furthermore, let $\FF'$ be the $\ZZ^{n}$-graded  free $S'$-resolution  of $K[H']$. Then, since $S$ is a flat $S'$-module,  $\FF'\tensor_{S'} S$ is a $\ZZ^n$-graded free $S$-resolution of $S/{I_{H'}}S$. The inclusion $S'/{I_{H'}}S\tensor_{S'}S\to S/{I_{H}}S$ induces a $\ZZ^n$-graded $S$-module complex homomorphism $\FF'\tensor_{S'}S\to \FF$. Applying $\_ \tensor_S K$ on this complex homomorphism with $K=S/\mm$, we obtain the following sequence of isomorphisms and  natural maps of $\ZZ^n$-graded $K$-modules

  \begin{multline*}
  \Tor_i^{S'}(K[H'],K)\iso H_i(\FF'\tensor_{S'}K)\iso H_i(\FF'\tensor_{S'}S)\tensor_SK)\to\\ 
 H_i(\FF\tensor_SK)\iso \Tor_i^S(K[H],K).
  \end{multline*}

 We need the following proposition several times in this paper. 
  
\begin{Proposition}\cite[Corollary\ 2.4]{[EHH15]}
  	\label{homologicallypure}
 Let $H'$ be a homologically pure subsemigroup  of $H$. If $\FF'$ is the minimal $\ZZ^n$-graded free $S'$-resolution of $K[H']$ and $\FF$ is the minimal $\ZZ^n$-graded free $S$-resolution of $K[H]$, then the complex homomorphism $\FF'\tensor S\to \FF$ induces an injective map $\FF'\tensor K\to \FF\tensor K$. Hence, $$\Tor_i^{S'}(K[H'],K)\to  \Tor_i^S(K[H],K)$$   is injective for all $i$. In particular, any minimal set of generators of $\Syz_i(K[H'])$ is part of a minimal set of generators  of  $\Syz_i(K[H])$. Moreover, $\beta_{ij}(K[H'])\leq \beta_{ij}(K[H])$ for all $i$ and $j$.
\end{Proposition}
  
Let us now define the semigroup ring structure on Hibi rings. Let $L=\MI(P)$ be a distributive lattice  with $P=\{p_1,\ldots,p_n\}$. For $ \alpha \in L$, define a $2n$-tuple $h_\alpha$ such that for $1 \leq i \leq n$, 
 \[
 \begin{cases} 
 1  & \text{at $i^{th}$ position if} \quad p_i \in \alpha  \text{,}  \\
 0 & \text{at $i^{th}$ position if} \quad p_i \notin \alpha \text{,}  \\
 0 & \text{at $(n+i)^{th}$ position if} \quad p_i \in \alpha \text{,}  \\
 1 & \text{at $(n+i)^{th}$ position if} \quad p_i \notin \alpha.
 \end{cases}
 \]
  
Let $H$ be the affine semigroup generated by $\{h_\alpha : \alpha \in L\}$. Then, we have $K[H]=R[L]$. Let $\beta, \gamma\in L$ such that $\beta \leq \gamma$. Define $L_1 = \{\alpha \in L : \beta \leq \ \alpha \leq\gamma\}$. Clearly, $L_1$ is a sublattice of $L$. Let $H_1$ be the affine subsemigroup of $H$ generated by $\{h_\alpha : \alpha \in L_1\}$.

\begin{Proposition}\label{semigroup}
	Let $H$ and $H_1$ be as defined above. Then $H_1$ is a  homologically pure subsemigroup of $H$.
\end{Proposition}
\begin{proof}
	We show that if $\alpha \notin L_1$ then $h - h_\alpha \notin H$ for all  $h \in H_1$. Suppose $\alpha \notin L_1$ then either $\alpha \nleq \gamma$ or $\alpha \ngeq \beta$.\par
	
	If $\alpha \nleq \gamma$, then there exists a $p_i \in \alpha$ such that $p_i \notin \gamma$. So $i^{th}$ entry of $h_\alpha$ is 1 but for any $\alpha' \in L_1$, $i^{th}$ entry of $h_{\alpha'}$ is 0. Hence, $h - h_\alpha \notin H$ for all  $h \in H_1$.\par
	
	If $\alpha \ngeq \beta$, then there exists a $p_j \in \beta$ such that $p_j \notin \alpha$. So $(n+j)^{th}$ entry of $h_{\alpha}$ is 1 but for any $\alpha' \in L_1$, $(n+j)^{th}$ entry of $h_{\alpha'}$ is 0. Hence, $h - h_\alpha \notin H$ for all  $h \in H_1$.
	
\end{proof}

\section{Property $N_p$ for Segre products of Hibi rings}\label{segre}

 In this section, we discuss the property $N_p$ of Segre products of Hibi rings for $p\in \{2,3\}$. Let $P_1$ and $P_2$ be two posets and $P$ be their disjoint union. From Section~\ref{hibiringdefi}, we know that the Segre product of Hibi rings $R[\MI(P_1)]$ and $R[\MI(P_2)]$ is isomorphic to the Hibi ring $R[\MI(P)]$. We also know that a  poset $\mathcal{P}$ is a chain if and only if $R[\MI(\mathcal{P})]$ is a polynomial ring. From Subsection~\ref{squarefree}, recall the definition of the semigroup associated to a Hibi ring. For $i\in \{1,2\}$, let $H_i$ be the affine semigroup generated by $\{h_\alpha : \alpha \in \MI(P_i)\}$ and let $H$ be the affine semigroup associated to the Hibi ring $R[\MI(P)]$. Since $\MI(P)=\{(\alpha,\beta): \alpha \in \MI(P_1)\ \text{and}\ \beta \in \MI(P_2) \}$, it is easy to see that $H$ is generated by $\{(h_{\alpha},h_{\beta}):  \alpha \in \MI(P_1)\ \text{and}\ \beta \in \MI(P_2)\}$.
 
\begin{Theorem}\label{unionbetti}
 Let $P_1,P_2,P, H_1,H_2$ and $H$ be as above. Then, for each $l\in \{1,2\}$, $H_l$ is isomorphic to a homologically pure subsemigroup of  $H$. In particular, if $\beta_{ij}(R[\MI(P_l)])\neq 0$ for some $l\in \{1,2\}$, then $\beta_{ij}(R[\MI(P)])\neq 0$.
\end{Theorem}
 
\begin{proof}
By symmetry, it suffices to prove the theorem for $l=1$. Consider the subsemigroup $G_1$ of $H$ generated by  $\{(h_{\alpha},h_{\emptyset}):  \alpha \in \MI(P_1)\}$, where $\emptyset$ is the minimal element of $\MI(P_2)$. It is easy to see that $G_1$ is isomorphic to the semigroup $H_1$. Also, observe that $\delta= (\emptyset,\emptyset)$ and $\gamma=(P_1,\emptyset)$  are the order ideals of $H$. The subsemigroup $G_1$ is generated by $\{h_\eta: \delta \leq \eta \leq \gamma\}$. So by Proposition~\ref{semigroup}, $G_1$ is a homologically pure subsemigroup of $H$. The second  part of the theorem follows from  Proposition~\ref{homologicallypure}. Hence the proof.
\end{proof} 
 
 \begin{Corollary}\label{segreposet}
 	Let $P$ be a poset such that it is a disjoint union of two posets $P_1$ and $P_2$.	If $R[\MI(P)]$ satisfies property $N_p$ for some $p\geq 2$, then so does $R[\MI(P_1)]$ and $R[\MI(P_2)]$. 
 \end{Corollary}
 
 \begin{proof}
 The proof  follows from  Theorem~\ref{unionbetti}. 
 \end{proof}

\begin{Lemma}\label{vanishingbetti}
Let $R[\MI(P)]$ be a Hibi ring associated to a poset $P$. Then the following statements hold:\\
$(a)$  If $\beta_{24}(R[\MI(P)])= 0$, then $R[\MI(P)]$ satisfies property $N_2$.\\
$(b)$  If $R[\MI(P)]$ satisfies property $N_2$ and $\beta_{35}(R[\MI(P)])= 0$, then it satisfies property $N_3$.
\end{Lemma} 
\begin{proof}
$(a)$ Hibi rings are algebra with straightening laws (ASL) and straightening relations are quadratic \cite[\S \ 2]{[HIBI87]}. ASL with quadratic straightening relations are Koszul \cite{[KEMPH90]}. So by \cite[Lemma\ 4]{[KEMPH90]}, $\beta_{2j}(R[\MI(P)])= 0$ for  all $j\geq 5$. This concludes the proof.\par
$(b)$ The proof follows from \cite[Theorem\ 6.1]{[ACI15]}.
\end{proof}

 \subsection{}\label{bettipoint}
 
In this  subsection, we  prove Theorem~\ref{beta24unionpoint}.
   
 \begin{Theorem}\label{beta24unionpoint}
 	Let $P_1$ be a poset,  $P_2 = \{b\}$ and $p\in \{2,3\}$. Let $P$ be the disjoint union of $P_1$ and $P_2$. If  $R[\MI(P_1)]$ satisfies property $N_p$ , then so does $R[\MI(P)]$.
 \end{Theorem}
 
  The proof of the above theorem follows the argument of Rubei \cite{[RUBEI02]}. Let $P_1$ and $P_2$ be as in theorem. So $\MI(P)=\{(\alpha,\beta): \alpha \in \MI(P_1),\ \beta \in \MI(P_2) \}$. Let $H$ be the affine semigroup generated by $\{(h_{\alpha},h_{\beta}):  \alpha \in \MI(P_1),\ \beta \in \MI(P_2)\}$. In order to prove the above theorem, by Proposition~\ref{bh} and Lemma~\ref{vanishingbetti}, it is enough to show that for $p\in \{2,3\}$, if $h =(h_1,h_2) \in H$ with $deg(h)=p+2$, then $\widetilde{H}_{p-1}(\Delta_h)= 0$. For $i=1,2$, let $H_i$ be the affine semigroup generated by $\{h_\alpha : \alpha \in \MI(P_i)\}$. Observe that $H_2$ is generated by two elements $h_{\emptyset}$ and $h_{\{b\}}$. For simplicity, we denote  $h_{\{b\}}$ by  $h_{b}$.
  
  Let $\Delta$ be a simplicial complex on a vertex set $V$. The {\em support} of a simplex $\sigma$ in $\Delta$ is the set of all vertices $v\in V$ such that $v\in \sigma$. Let $\alpha=\sum_{i}^{}a_i\sigma_i$ where $c_i\in \ZZ$, be a chain in $\Delta$. The {\em support} of $\alpha$, denoted by $sp(\alpha)$, is the union of the support of the simplexes $\sigma_i$. We  denote the {\em $i$-skeleton} of $\Delta$ by $sk^i(\Delta)$.
 
 \begin{Notation}\label{f_delta_g}
 	 Let $g\in H_1$ with $deg(g) = d$.\\
 	$(a)$ Denote $g_\epsilon = (g,g')$, where $g' = (d-\epsilon) h_{\emptyset}+\epsilon h_{b}\in H_2$ and $\epsilon\in \{0,\ldots,d\}$.\\
 	$(b)$ For  $0\leq l\leq d-1$, let
    $$F^l(\Delta_g)= \mathop{\cup}_{\substack{g_1,\ldots, g_d \\ \text{s.t.}\ g_1+\ldots+g_d=g}}^{}\mathop{\cup}_{i_0,...,i_l\in\{1,\ldots,d\}}^{}\big\langle (g_{i_{0}},h_{\emptyset}),\ldots, (g_{i_{l}},h_{\emptyset})\big \rangle .$$ 
 \end{Notation}
 \begin{Lemma}\label{homology}
 	Under the notations of Notation~\ref{f_delta_g}.\\ 
 	$(a)$ For all $i\leq l-1$,  $\widetilde{H}_{i}(F^l(\Delta_g)) \cong \widetilde{H}_{i}(\Delta_g)$.\\
 	$(b)$ For $\epsilon \in \{1,2\}$, $F^l(\Delta_g) \subseteq \Delta_{g_{\epsilon}}$ if and only if $l\leq d-\epsilon -1$.
 \end{Lemma}
 \begin{proof}
 $(a)$	The proof follows from $F^l(\Delta_g) \cong sk^l(\Delta_g)$.\\
 $(b)$  Let $g_1,\ldots,g_d\in H_1$ be such that $\sum_{i=1}^{d} g_i=g$.  Observe that for any $\{i_0,...,i_l\}\subseteq[d]$, \ab $\{(g_{i_{0}},h_{\emptyset}),\ldots, (g_{i_{l}},h_{\emptyset})\}$ is a simplex in $\Delta_{g_{\epsilon}}$ if and only if  $l\leq d-\epsilon -1$.
 \end{proof}
 
 Let $g\in H_1$ with $deg(g)=d$ and let $\epsilon\in\{0,\ldots,d\}$.  Note that $\Delta_{g_{\epsilon}} \cong \Delta_{g}$ for all $\epsilon\in \{0,d\}$. Also, we have $\Delta_{g_{\epsilon}} \cong \Delta_{g_{d-\epsilon}}$. Thus, to  prove the theorem, it suffices to consider the cases $h_2 = (p+2-\epsilon) h_{\emptyset}+\epsilon h_{b}$, where $\epsilon\in\{1,2\}$ and $p\in \{2,3\}$.
 
 \begin{Remark}\label{cone}
 	Let $g\in H_1$ with $deg(g)=d$ and $\epsilon\in \{0,\ldots,d\}$. Let $g_1,...,g_d\in H_1$ be such that $g = \sum_{i=1}^{d}g_i$. Observe that $\sigma =\big\{ (g_{i_{1}},h_{\emptyset}),\ldots, (g_{i_{d-\epsilon+1}} ,h_{\emptyset})\big\} \notin \Delta_{g_{\epsilon}}$ for any $i_1,\ldots, i_{d-\epsilon+1}\in \{1,\ldots,d\}$. For $l \in \{1,\ldots,d\}$ with $l\neq i_j$, $j \in \{1,\ldots,d-\epsilon+1\}$, let 
 	$$\sigma' = \sum_{j=1}^{d-\epsilon+1}(-1)^{j-1} \big\{ (g_{i_{l}},h_{b}),(g_{i_{1}} ,h_{\emptyset}),\ldots, \widehat{(g_{i_{j}} ,h_{\emptyset})},\ldots,(g_{i_{d-\epsilon+1}} ,h_{\emptyset})\big\}$$
 	be a $(d-\epsilon)$-chain in $\Delta_{g_{\epsilon}}$. Then $\partial \sigma = \partial \sigma'$.
 \end{Remark}

 \begin{Definition}
 	For any $g \in H_1$ with $deg(g) = d$ and $\epsilon \in \{1,\ldots, d\}$, we define $R_{g,\epsilon}$ to be  the following simplicial complex:
    $$\mathop{\cup}_{\substack{g_1,\ldots, g_d\in H_1 \ \\ \text{s.t.}\ g_1+\ldots+g_d=g}}^{}\mathop{\cup}_{\substack{i_1,...,i_{d-1}\in\{1,\ldots,d\}\\
    i_l\neq i_m}}^{}\big \langle (g_{i_{1}},h_{b}),\ldots, (g_{i_{\epsilon-1}} ,h_{b}),(g_{i_{\epsilon}},h_{\emptyset}),\ldots,(g_{i_{d-1}},h_{\emptyset})\big\rangle .$$ 
 \end{Definition}
 
 \begin{Lemma}\label{rgepsilon}
 	Let $g\in H_1$ with $deg(g) = d$ and $\epsilon \in \{1,2\}$.  Assume that\par
 	$(a) \  (i,d) \in \{(0,3),(1,4)\}$ and\par
 	$(b) \ \widetilde{H}_{i}(\Delta_{g_{\epsilon-1}})= 0$.\\
 Then	$\widetilde{H}_{i}(R_{g,\epsilon})= 0$.
 \end{Lemma}
 
 \begin{proof}
 	Observe that $R_{g,\epsilon} \subseteq \Delta_{g_{\epsilon-1}}$. If $\epsilon=1$, then $sk^2(\Delta_{g_{\epsilon-1}})\subseteq sk^2(R_{g,\epsilon})$. Thus, $\widetilde{H}_{i}(\Delta_{g_{\epsilon-1}}) = \widetilde{H}_{i}(R_{g,\epsilon})$ for $i=0,1$. So we only have to consider the case $\epsilon=2$. Let $\beta$ be an $i$-cycle in $R_{g,\epsilon}$. Since $\widetilde{H}_{i}(\Delta_{g_{\epsilon-1}})= 0$, there exists an $(i+1)$-chain $\eta$ in $\Delta_{g_{\epsilon-1}}$ such that $\partial\eta =\beta$. Suppose $\eta = \sum_{j}^{}c_j\sigma_j$, where $\sigma_j$ is an $(i+1)$-simplex in $\Delta_{g_{\epsilon-1}}$. Now consider an $(i+1)$-chain $\psi$ in $R_{g,\epsilon}$ such that $\psi = \sum_{j}^{}c_j\sigma'_j$, where $\sigma'_j = \sigma_j$ if $\sigma_j \in  R_{g,\epsilon}$ else $\sigma'_j$ is as defined in Remark~\ref{cone} corresponding to $\sigma_j$. Then $\partial\psi = \beta$.
 \end{proof}	
 
 \begin{Lemma}\label{1cyclehomo}
 	Let $g \in H_1$ with $deg(g)=4$ and $\epsilon \in \{1,2\}$. Every 1-cycle $\gamma$ in  $\Delta_{g_{\epsilon}}$ is  homologous to  an 1-cycle in $F^1(\Delta_g)\ (\subseteq \Delta_{g_{\epsilon}})$.	
 \end{Lemma}
 \begin{proof}
 	We prove the lemma by induction on the cardinality of $(sp(\gamma)\cap sk^0(\Delta_{g_{\epsilon}}))\setminus F^1(\Delta_g).$
 	
 	Let $(f,h_b) \in sp(\gamma)$. Let $\MS_{(f,h_b)}$ be the set of $1$-simplexes of $\gamma$ with vertex $(f,h_b)$. For $\sigma=\{v,(f,h_b)\} \in \MS_{(f,h_b)}$, let $\sigma'= \{v,(f,h_{\emptyset})\}$. Clearly, $\sigma'$ is an 1-simplex of $\Delta_{g_{\epsilon}}$. Let $\alpha = \sum_{\sigma \in \MS_{(f,h_b)}}^{}(-\sigma+\sigma')$ be the 1-cycle in $\Delta_{g_{\epsilon}}$.\par 
 	
 	Now, we show that for a vertex $v$  in $\Delta_{g_{\epsilon}}$, if $\langle v,(f,h_b)\rangle \subseteq \Delta_{g_{\epsilon}}$, then $v\in R_{g-f,\epsilon}$.
 	Observe that  
    $$R_{g-f,1} =\mathop{\cup}_{\substack{g_1,g_2, g_3\in H_1 \\ \text{s.t.}\ g_1+g_2+g_3=g-f}}^{}\mathop{\cup}_{\substack{i_1,i_{2}\in\{1,2,3\}\\
    i_1\neq i_2}}^{}\big \langle(g_{i_{1}},h_{\emptyset}),(g_{i_{2}},h_{\emptyset})\big\rangle $$
 	and  
    $$R_{g-f,2} =\mathop{\cup}_{\substack{g_1,g_2, g_3\in H_1 \\ \text{s.t.}\ g_1+g_2+g_3=g-f}}^{}\mathop{\cup}_{\substack{i_1,i_{2}\in\{1,2,3\}\\
    i_1\neq i_2}}^{}\big \langle(g_{i_{1}},h_{b}),(g_{i_{2}},h_{\emptyset})\big\rangle.$$
 	
 	If $\epsilon =1$, then $v = (f',h_{\emptyset})$ for some $f'\in H_1$ such that $g-(f+f')\in H_1$. If $\epsilon =2$, then either $v = (f',h_{\emptyset})$ or $v = (f',h_b)$ for some $f'\in H_1$ such that $g-(f+f')\in H_1$. In both cases, $v\in R_{g-f,\epsilon}$.\par
 	So we obtain that $sp(\alpha)\subseteq C$, where $C$ is the union of the cones $\langle (f,h_b),R_{g-f,\epsilon} \rangle$ and $\langle (f,h_{\emptyset}),R_{g-f,\epsilon} \rangle$. Notice that $C\subseteq \Delta_{g_{\epsilon}}$. Since Hibi rings satisfy property $N_1$, we have $\widetilde{H}_{0}(\Delta_{({g-f)}_{\epsilon-1}})= 0$. Thus, $\widetilde{H}_{0}(R_{g-f,\epsilon})= 0$   by Lemma~\ref{rgepsilon}. Since $\widetilde{H}_{i}(C)=\widetilde{H}_{i-1}(R_{g-f,\epsilon})$, we have $\widetilde{H}_{1}(C)=0$. Thus, $\alpha$ is homologous to 0 which implies that $\sigma$ is homologous to $\sigma+\alpha$.\par
 	
 	Observe that $\#((sp(\gamma+\alpha)\cap sk^0(\Delta_{g_{\epsilon}}))\setminus F^1(\Delta_g))<\#((sp(\gamma)\cap sk^0(\Delta_{g_{\epsilon}}))\setminus F^1(\Delta_g)).$ Hence, we conclude the proof by induction.
 \end{proof}
 
 \begin{proof}[Proof of Theorem~\ref{beta24unionpoint} for $N_2$]
 	We have to show that if $h =(h_1,h_2) \in H$  with $deg(h)=4$, then $\widetilde{H}_{1}(\Delta_h)= 0$. We consider the following cases:
 	\begin{itemize}
 		\item[$(1).$] Assume that $h_2 = 3 h_{\emptyset}+h_{b}$. Let $\gamma$ be an 1-cycle in $\Delta_h$. By Lemma~\ref{1cyclehomo}, $\gamma$ is  homologous to  an 1-cycle $\gamma'$ of $F^1(\Delta_{h_{1}})\subset\Delta_{h}$. In other words, there exists a 2-chain $\mu$ in $\Delta_h$ such that $\partial\mu = \gamma-\gamma'$. Also, $\widetilde{H}_{1}(F^2(\Delta_{h_{1}})) \cong \widetilde{H}_{1}(\Delta_{h_{1}})=0$, where the isomorphism is due to Lemma~\ref{homology}(a) and the equality is by hypothesis. As $F^1(\Delta_{h_{1}})\subset F^2(\Delta_{h_{1}})$, there exists a 2-chain $\mu'$ in $F^2(\Delta_{h_{1}})$ such that $\partial\mu' = \gamma'$. Since $F^2(\Delta_{h_{1}})\subseteq \Delta_h$, $\mu'$ is a 2-chain in $\Delta_h$. Therefore, $[\gamma'] = 0$ in $\widetilde{H}_{1}(\Delta_{h})$. So $[\gamma] = 0$ in $\widetilde{H}_{1}(\Delta_{h})$.   
 		
 		\item[$(2).$] Consider $h_2 = 2 h_{\emptyset}+2h_{b}$. Every 1-cycle $\gamma$ in  $\Delta_h$ is  homologous to  an 1-cycle $\gamma'$ in $F^1(\Delta_{h_{1}})$ by Lemma~\ref{1cyclehomo}. But in this case, $F^2(\Delta_{h_{1}})$ is not contained in $\Delta_h$. Since $\widetilde{H}_{1}(F^2(\Delta_{h_{1}})) = 0$, there exists a  2-chain $\mu$  in $F^2(\Delta_{h_{1}})$ such that $\partial\mu = \gamma'$. Let $\mu = \sum_{i}^{}c_i\sigma_i$, where $\sigma_i$ is a 2-simplex in $F^2(\Delta_{h_{1}})$. Consider a 2-chain $\psi$ in $\Delta_h$ such that $\psi = \sum_{i}^{}c_i\sigma'_i$, where $\sigma'_i =\sigma_i$ if $\sigma_i \in \Delta_h$ else $\sigma'_i$ is as defined in Remark~\ref{cone} corresponding to $\sigma_i$. Then $\partial\psi = \gamma'$. Therefore, $[\gamma']=0$ in $\widetilde{H}_{1}(\Delta_h)$.  So $[\gamma]=0$ in $\widetilde{H}_{1}(\Delta_h).$ Hence the proof.
 	\end{itemize}
 \end{proof}
 
 \begin{Lemma}\label{2cyclehomo}
 	Let $g \in H_1$ with $deg(g)=5$ and $\epsilon \in \{1,2\}$. Every 2-cycle $\gamma$ in  $\Delta_{g_{\epsilon}}$ is  homologous to  a 2-cycle in $F^2(\Delta_g)  \ (\subset\Delta_{g_{\epsilon}})$.
 \end{Lemma}
 \begin{proof}
 	We prove the result by induction on the cardinality of $(sp(\gamma)\cap sk^0(\Delta_{g_{\epsilon}}))\setminus F^2(\Delta_g).$\par
 	
 	Let $(f,h_b) \in sp(\gamma)$. Let $\MS_{(f,h_b)}$ be the set of $2$-simplexes of $\gamma$ with vertex $(f,h_b)$. For $\sigma=\{v,u,(f,h_b)\} \in \MS_{(f,h_b)}$, let $\sigma'= \{v,u,(f,h_{\emptyset})\}$. Clearly, $\sigma'$ is a 2-simplex of $\Delta_{g_{\epsilon}}$. Let $\alpha = \sum_{\sigma \in \MS_{(f,h_b)}}^{}(-\sigma+\sigma')$ be the 2-cycle in $\Delta_{g_{\epsilon}}$.\par 
 	Now, we show that for vertexes $v,u$  in $\Delta_{g_{\epsilon}}$, if $\langle v,u,(f,h_b)\rangle \subseteq \Delta_{g_{\epsilon}}$, then $\langle v,u \rangle \subseteq R_{g-f,\epsilon}$.
 	Observe that  
    $$R_{g-f,1} =\mathop{\cup}_{\substack{g_1,\ldots, g_4\in H_1 \\ \text{s.t.}\ g_1+\ldots+g_4=g-f}}^{}\mathop{\cup}_{\substack{i_1,i_{2},i_{3}\in\{1,\ldots,4\}\\
    i_l\neq i_k}}^{}\big \langle(g_{i_{1}},h_{\emptyset}),(g_{i_{2}},h_{\emptyset})(g_{i_{3}},h_{\emptyset})\big\rangle$$
 	and  
    $$R_{g-f,2} =\mathop{\cup}_{\substack{g_1,\ldots, g_4\in H_1 \\ \text{s.t.}\ g_1+\ldots+g_4=g-f}}^{}\mathop{\cup}_{\substack{i_1,i_{2},i_{3}\in\{1,\ldots,4\}\\
    i_l\neq i_k}}^{}\big \langle(g_{i_{1}},h_{b}),(g_{i_{2}},h_{\emptyset})(g_{i_{3}},h_{\emptyset})\big\rangle .$$
 	
 	If $\epsilon =1$, then $\langle v,u \rangle = \langle (f_1,h_{\emptyset}),(f_2,h_{\emptyset}) \rangle$ for some $f_1,f_2\in H_1$ such that $g-(f+f_1+f_2)\in H_1$. If $\epsilon =2$, then either $\langle v,u \rangle = \langle (f_1,h_{\emptyset}),(f_2,h_{\emptyset}) \rangle$ or $\langle v,u \rangle = \langle (f_1,h_b),(f_2,h_{\emptyset}) \rangle$ for some $f_1,f_2\in H_1$ such that $g-(f+f_1+f_2)\in H_1$. In both cases, $\langle v,u \rangle = \langle (f_1,h_{\emptyset}),(f_2,h_{\emptyset}) \rangle  \subseteq R_{g-f,\epsilon}$.\par
 	
 	So we obtain that $sp(\alpha)\subseteq C$, where $C$ is the union of the cones $\langle (f,h_b),R_{g-f,\epsilon} \rangle$ and $\langle (f,h_{\emptyset}),R_{g-f,\epsilon} \rangle$. Notice that $C\subseteq \Delta_{g_{\epsilon}}$. Since $R[\MI(P)]$ satisfies property $N_2$, we have $\widetilde{H}_{1}(\Delta_{({g-f)}_{\epsilon-1}})= 0$. Thus by Lemma~\ref{rgepsilon}, $\widetilde{H}_{1}(R_{g-f,\epsilon})= 0$. Since $\widetilde{H}_{i}(C)=\widetilde{H}_{i-1}(R_{g-f,\epsilon})$, we have $\widetilde{H}_{2}(C)=0$. Thus, $\alpha$ is homologous to 0. So $\sigma$ is homologous to $\sigma+\alpha$.\par
 	
 	Observe that $\#((sp(\gamma+\alpha)\cap sk^0(\Delta_{g_{\epsilon}}))\setminus F^2(\Delta_g))<\#((sp(\gamma)\cap sk^0(\Delta_{g_{\epsilon}}))\setminus F^2(\Delta_g)).$ Hence, we conclude the proof by induction.
 \end{proof}
 
 \begin{proof}[Proof of Theorem~\ref{beta24unionpoint} for $N_3$]
 	We have to show that if $h =(h_1,h_2) \in H$ with $deg(h)=5$, then $\widetilde{H}_{2}(\Delta_h)= 0$. We consider the following cases:
 	\begin{itemize}
 		\item[$(1).$] Consider $h_2 = 4 h_{\emptyset}+h_{b}$. Let $\gamma$ be a 2-cycle in $\Delta_h$. By Lemma~\ref{2cyclehomo}, $\gamma$ is  homologous to  a 2-cycle $\gamma'$ of $F^3(\Delta_{h_{1}})\subset\Delta_{h}$. In other words, there exists a 3-chain $\mu$ in $\Delta_h$ such that $\partial\mu = \gamma-\gamma'$. Also, $\widetilde{H}_{2}(F^3(\Delta_{h_{1}})) \cong \widetilde{H}_{2}(\Delta_{h_{1}})=0$, where the isomorphism is due to Lemma~\ref{homology} and the equality holds because $R[\MI(P_1)]$ satisfies property $N_3$. As $F^2(\Delta_{h_{1}})\subset F^3(\Delta_{h_{1}})$, there exists a 3-chain $\mu'$ in $F^3(\Delta_{h_{1}})$ such that $\partial\mu' = \gamma'$. Since $F^3(\Delta_{h_{1}})\subseteq \Delta_h$, $\mu'$ is a 3-chain in $\Delta_h$. Therefore, $[\gamma'] = 0$ in $\widetilde{H}_{1}(\Delta_{h})$. So $[\gamma] = 0$ in $\widetilde{H}_{1}(\Delta_{h})$.   
 		
 		\item[$(2).$] Assume $h_2 = 3h_{\emptyset}+2h_{b}$. By Lemma~\ref{2cyclehomo}, every 2-cycle $\gamma$ in  $\Delta_h$ is  homologous to  a 2-cycle $\gamma'$ in $F^2(\Delta_{h_{1}})$. But in this case, $F^3(\Delta_{h_{1}})\nsubseteq \Delta_h$. Since $\widetilde{H}_{2}(F^3(\Delta_{h_{1}})) = 0$, there exists a  3-chain $\mu$  in $F^3(\Delta_{h_{1}})$ such that $\partial\mu = \gamma'$. Let $\mu = \sum_{i}^{}c_i\sigma_i$, where $\sigma_i$ is a 3-simplex in $F^3(\Delta_{h_{1}})$. Consider a 3-chain $\psi$ in $\Delta_h$ such that $\psi = \sum_{i}^{}c_i\sigma'_i$, where $\sigma'_i =\sigma_i$ if $\sigma_i \in \Delta_h$ else $\sigma'_i$ is as defined in Remark~\ref{cone} corresponding to $\sigma_i$. Then $\partial\psi = \gamma'$. Therefore, $[\gamma']=0$ in $\widetilde{H}_{2}(\Delta_h)$. So $[\gamma]=0$ in $\widetilde{H}_{2}(\Delta_h).$ Hence the proof.
 	\end{itemize}
 \end{proof}
 
\subsection{}\label{bettichain}
In this subsection, we show that if a Hibi ring satisfies property $N_2$, then its Segre product with a polynomial ring in finitely many variables also satisfies property $N_2$.

\begin{Proposition}\label{bettiunion}
Let $P$ be a poset such that it is a disjoint union of a poset $P_1$ and a chain $P_2=\{a_1,\ldots,a_n\}$ with $a_1\lessdot\cdots\lessdot a_n$. Let $\{x\}$ be a poset and $P_2'$ be the ordinal sum $P_2\dirsum\{x\}$. Let $Q$ be the disjoint union of the posets $P_1$ and $P_2'$.	If $R[\MI(P)]$ satisfies property $N_2$, then so does $R[\MI(Q)]$.    \end{Proposition}
   
 \begin{Theorem}\label{chain}
 Let $P$ be a poset such that it is a disjoint union of a poset $P_1$ and a chain $P_2$. If $R[\MI(P_1)]$ satisfies property $N_2$, then so does $R[\MI(P)]$. 
 \end{Theorem}
  
 \begin{proof}
 The proof follows from Theorem~\ref{beta24unionpoint} and Proposition~\ref{bettiunion}.
 \end{proof}

 Now, this subsection is dedicated to the proof of the Proposition~\ref{bettiunion}. The proof of Proposition~\ref{bettiunion} is motivated from Rubei \cite{[RUBEI07]}.\par
  
 Let $P$ be a poset such that it is a disjoint union of two posets $P_1$ and $P_2$. Let $\{x\}$ be a poset and $P_2'$ be the ordinal sum $P_2\dirsum\{x\}$. Let $Q$ be the disjoint union of posets $P_1$ and $P_2'$. Let $H$ and $H'$  be the affine semigroups generated by $\{h_{\alpha}:  \alpha \in \MI(Q)\}$ and $\{h_{\beta}:  \beta \in \MI(P)\}$ respectively. For $i\in \{1,2\}$, let $H_i$ be the affine semigroup generated by $\{h_\alpha : \alpha \in \MI(P_i)\}$. For $\alpha\in Q$, the first entry of $h_\alpha$ is 1 if $x\in \alpha$ and the second entry of $h_{\alpha}$ is 1 if $x\notin \alpha$. 
 
 $Note:$ Let  $h\in H$ with $deg(h)=d$. In this subsection, we either denote $h$ by $(\epsilon, d-\epsilon, h')$, where $h'\in H'$ or we denote it by $(\epsilon, d-\epsilon, h_2, h_1)$, where $h_i\in H_i$ for all $i=1,2$.
  
 Let $h\in H$ with $deg(h)=d$. Then $h = (\epsilon, d-\epsilon, h')$, where $h'\in H'$, $\epsilon\leq d,$ $\epsilon\in \NN$. Let $X_h$ be the  following simplicial  complex:
  $$X_h:=\Delta_h \cup \Delta_{(\epsilon-1,d-\epsilon+1,h')}\cup \ldots \cup \Delta_{(0,d,h')}.$$
  Observe that $\Delta_{(0,d,h')} \cong \Delta_{h'}\cong \Delta_{(d,0,h')}$. Let $\gamma$ be an $1$-cycle in $X_h$. For every vertex $v\in \gamma$, let $\mathcal{S}_{v,\gamma}$ be the set of simplexes of $\gamma$ with vertex $v$ and $\mu_{v,\gamma}$ be the $0$-cycle such that $v*\mu_{v,\gamma} = \sum_{\tau\in \mathcal{S}_{v,\gamma}}^{}\tau$, where $*$ denotes the joining. 
  
  \begin{Proposition}\label{homologous}
  	Let $h\in H$ with $deg(h)=4$. Let $\gamma$ be  an 1-cycle in $\Delta_h$. Then there exists an 1-cycle $\gamma'$ in $\Delta_{(0,4,h')}$ such that $\gamma$ is homologous to $\gamma'$ in $X_h$. 	
  \end{Proposition}
  
  \begin{proof}
  	Let $h_{\alpha_{1}},\ldots, h_{\alpha_{m}}$ be the vertices of $\gamma$ with non-zero first entry. In other words, these are all the vertices $h_{\alpha}$ of $\gamma$ such that $x\in \alpha$. For $1\leq i\leq m$, let $\beta_{i} = \alpha_i \setminus \{x\}$. 
  	Observe that $\mu_{h_{\alpha_{1}},\gamma}$ is in $\Delta_{h-h_{\alpha_{1}}}$ and 
  	$$\widetilde{H}_{1}(h_{\alpha_{1}} * \Delta_{h-h_{\alpha_{1}}}\cup h_{\beta_{1}} * \Delta_{h-h_{\alpha_{1}}})= \widetilde{H}_{0}(\Delta_{h-h_{\alpha_{1}}})= 0,$$
  	where the last equality holds because $R[\MI(Q)]$ satisfies property $N_1$. So $h_{\alpha_{1}} * \mu_{h_{\alpha_{1}},\gamma}- h_{\beta_{1}} * \mu_{h_{\alpha_{1}},\gamma}$ is homologous to 0 in $X_h$. Hence, $\gamma_1:= \gamma -(h_{\alpha_{1}} * \mu_{h_{\alpha_{1}},\gamma}- h_{\beta_{1}} * \mu_{h_{\alpha_{1}},\gamma})$ is homologous to $\gamma$ in $X_h$. Informally speaking, we have got $\gamma_1$ from $\gamma$ by replacing the vertex $h_{\alpha_{1}}$ with $h_{\beta_{1}}$.
  	Inductively, define
  	$$\gamma_i := \gamma_{i-1}-(h_{\alpha_{i}} * \mu_{h_{\alpha_{i}},\gamma_{i-1}}- h_{\beta_{i}} * \mu_{h_{\alpha_{i}},\gamma_{i-1}})$$
  	for $2\leq i \leq m$. Since all the vertexes of $\gamma_m$ have first entry zero, we have $\gamma_m \in \Delta_{(0,d,h')}$. We set $\gamma' =\gamma_m$ and prove that $\gamma_m$ is homologous to $\gamma$ in $X_h$. To prove this, it suffices to show that $h_{\alpha_{i}} * \mu_{h_{\alpha_{i}},\gamma_{i-1}}- h_{\beta_{i}} * \mu_{h_{\alpha_{i}},\gamma_{i-1}}$  is homologous to 0 for $2\leq i \leq m$.
  	
  	Let $\theta_0$ be the sum of simplexes $\tau$ of $\mu_{h_{\alpha_{i}},\gamma_{i-1}}$ such that  $\tau$ is a vertex of $\Delta_{h-h_{\alpha_{i}}}$ and let $\theta_1$ be the sum of simplexes $\tau$ of $\mu_{h_{\alpha_{i}},\gamma_{i-1}}$ such that  $\tau$ is not a vertex of $\Delta_{h-h_{\alpha_{i}}}$. Observe that $\mu_{h_{\alpha_{i}},\gamma_{i-1}}= \theta_0 +\theta_1$ and  $\theta_1$ is a 0-cell in $\Delta_{(\epsilon-1,5-\epsilon,h')-h_{\alpha_{i}}}$. Since $\mu_{h_{\alpha_{i}},\gamma}$ is a 0-cycle, $\theta_0$ is a 0-cycle of $\Delta_{h-h_{\alpha_{i}}}$ and $\theta_1$ is a 0-cycle of $\Delta_{(\epsilon-1,5-\epsilon,h')-h_{\alpha_{i}}}$. Furthermore, since $R[\MI(Q)]$ satisfies property $N_1$, $\theta_0$ is homologous to 0 in $\Delta_{h-h_{\alpha_{i}}}$ and $\theta_1$ is homologous to 0 in $\Delta_{(\epsilon-1,5-\epsilon,h')-h_{\alpha_{i}}}$. Thus, they are homologous to 0 in $X_h$.  Therefore, $h_{\alpha_{i}} * \mu_{h_{\alpha_{i}},\gamma_{i-1}}- h_{\beta_{i}} * \mu_{h_{\alpha_{i}},\gamma_{i-1}}$  is homologous to 0 in $X_h$. This concludes the proof.
  \end{proof}
  
  \begin{Lemma}\label{chaincomplex}
  	Let $\mathcal{P}=\{q_1,\ldots,q_r\}$  be a chain  such that $q_1\lessdot\cdots\lessdot q_r$ and let $\mathcal{H}$ be the semigroup corresponding to $R[\MI(\mathcal{P})]$. Let $h = \sum_{i=1}^{d}h_{\alpha_{i}}\in \mathcal{H}$. Then 
  	$$\Delta_h =\langle h_{\alpha_{1}},\ldots, h_{\alpha_{d}} \rangle.$$
  \end{Lemma}
  
  \begin{proof}
  	It is enough to show that for some $\alpha\in \MI(\mathcal{P})$, if $h_{\alpha} \notin \{h_{\alpha_{1}},\ldots, h_{\alpha_{d}}\}$, then $h_{\alpha}$ is not a vertex of $\Delta_h$. If $\alpha=\emptyset$, then the entry corresponding to ``$q_1\notin\alpha$'' in $h-h_\alpha$ will be -1. So $h_{\alpha}$ is not a vertex of $\Delta_h$. If $\alpha_i\leq \alpha$ for all $i\in \{1,\ldots,d\}$, then $h-h_{\alpha}\notin \mathcal{H}$. Hence, $h_{\alpha}$ is not a vertex of $\Delta_h$. Now suppose for all $i\in \{1,\ldots,d\}$, $\alpha_i \nleq \alpha$. Let $\{h_{\alpha_{i_{1}}},\ldots, h_{\alpha_{i_{m}}}\}$ be the subset of $\{h_{\alpha_{1}},\ldots, h_{\alpha_{d}}\}$ such that $ h_{\alpha}< h_{\alpha_{i_{j}}}$ for all $j=1,\ldots,m$. Let $\alpha = \{q_1,\ldots,q_s\}$, where $1\leq s\leq  r-1$. Observe that the entries corresponding to $q_s$  and $q_{s+1}$ in $h$ are $m$. But the entries corresponding $q_s$ and $q_{s+1}$ in $h-h_\alpha$ are $m-1$ and $m$ respectively. Hence, $h-h_\alpha\notin \mathcal{H}$.  This completes the proof. 
  \end{proof}
  
  From now onwards, let $P$ be a poset such that it is a disjoint union of a poset $P_1$ and a chain $P_2=\{a_1,\ldots,a_n\}$ with $a_1\lessdot\cdots\lessdot a_n$. Let $P_2'$ be the ordinal sum $P_2\dirsum\{x\}$. Furthermore, let $Q$ be the disjoint union of posets $P_1$ and $P_2'$. Let $H_{P'_{2}}$ be the  semigroup associated to $R[\MI(P'_2)]$.

  \begin{Lemma}\label{counting}
  	Let $h = \sum_{i=1}^{d}h_{\alpha_{i}}\in H$ with $h =(\epsilon,d-\epsilon,h')$, $\epsilon\geq 1$. For  $r<d$, assume that there are  exactly $r$ number of $i$'s with $\alpha_i= \alpha_i^1\cup P_2$, where $\alpha_i^1\in \MI(P_1)$. Let $\tau =\{h_{\beta_{1}},\ldots,h_{\beta_{m+1}}\}$ be an $m$-simplex of $X_h$. Then $\tau \in  \Delta_h$ if and only if there are atmost $r$ number of $\beta_i$'s with $\beta_i= \beta_i^1\cup P_2$, where $\beta_i^1\in \MI(P_1)$. 
  \end{Lemma}
  
  \begin{proof}
  	If we write $h = (\epsilon,d-\epsilon,h_2, h_1)$, where $h_i\in H_{i}$ for all $i=1,2$, then $(\epsilon-1,d-\epsilon+1,h') = (\epsilon-1,d-\epsilon+1,h_2, h_1)$. For $1\leq i\leq d$, let $\alpha_{i}=(\alpha_{i}^2,\alpha_{i}^1)$, where $\alpha_{i}^2\in \MI(P_2')$ and $\alpha_{i}^1\in \MI(P_1)$, So we can write $h_{\alpha_{i}} = (1,0,h_{P_{2}},h_{\alpha_{i}^1})$ if $x\in \alpha_{i}$ and $h_{\alpha_{i}} = (0,1,h_{\alpha_{i}^2},h_{\alpha_{i}^1})$ if $x\notin \alpha_{i}^2$, where $h_{\alpha_{i}^2}\in H_{2}$ and $h_{\alpha_{i}^1}\in H_{1}$. 
  	We have 
  	$$h_1 = \sum_{i=1}^{d}h_{\alpha^1_{i}},\quad (\epsilon,d-\epsilon,h_2)=\sum_{i=1}^{d}h_{\alpha^2_{i}}.$$ Let 
  	$\tau_1=\{h_{\beta_{1}^1},\ldots,h_{\beta_{m+1}^1}\}$ and $\tau_2=\{h_{\beta_{1}^2},\ldots,h_{\beta_{m+1}^2}\}$. Note that $\tau_1$, $\tau_2$ could be multisets.\par
  	Now we show that  $\tau\in \Delta_h$ if and only if  $\tau_2 \subseteq \{h_{\alpha_{1}^2},\ldots,h_{\alpha_{d}^2}\}$. 	Observe that if $\tau\in \Delta_h$, then $(\epsilon,d-\epsilon,h_2)- \sum_{j=1}^{m+1}h_{\beta_{j}^2}\in H_{P'_{2}}$. Hence, $\tau_2 \subseteq \{h_{\alpha_{1}^2},\ldots,h_{\alpha_{d}^2}\}$, by Lemma~\ref{chaincomplex}. On the other hand, if $\tau_2 \subseteq \{h_{\alpha_{1}^2},\ldots,h_{\alpha_{d}^2}\}$, then $(\epsilon,d-\epsilon,h_2)- \sum_{j=1}^{m+1}h_{\beta_{j}^2}\in H_{P'_{2}}$. Since $\tau\in X_h$, there exists an $i_0\in \{0,\ldots,\epsilon\}$  such that $\tau\in \Delta_{(\epsilon-i_0,d-\epsilon+i_0,h_2,h_1)}.$ So  $(\epsilon-i_0,d-\epsilon+i_0,h_2,h_1)- \sum_{j=1}^{m+1}h_{\beta_{j}}\in H$. Therefore,  $h_1- \sum_{j=1}^{m+1}h_{\beta_{j}^1}\in H_{1}$. We obtain $\tau \in \Delta_h$.\par
  	
  	The proof of `only if' part follows from the above claim. To prove `if', it suffices to show that $\tau_2 \subseteq \{h_{\alpha_{1}^2},\ldots,h_{\alpha_{d}^2}\}$. For $1\leq i \leq \epsilon$,  we have	
  	$$(\epsilon-i,d-\epsilon+i,h_2)=\sum_{j=1}^{\epsilon-i}h_{P'_2}+ \sum_{j=\epsilon-i+1}^{\epsilon}h_{P_2}+\sum_{j=\epsilon+1}^{\epsilon+r}h_{P_2}+\sum_{j=\epsilon-r+1}^{d}h_{\alpha^2_{j}}\in H_{P'_{2}}.$$
  	Let $i_0\in \{0,\ldots,\epsilon\}$ be such that $\tau\in \Delta_{(\epsilon-i_0,d-\epsilon+i_0,h_2,h_1)}.$ Since there are atmost $r$ number of $\beta_i$'s in $\tau$ with $\beta_i= \beta_i^1\cup P_2$, where $\beta_i^1\in \MI(P_1)$, the multiplicity of $h_{P_{2}}$ in $\tau_2$ is atmost $r$. So by Lemma~\ref{chaincomplex}, 
  	$$\tau_2\subseteq \{h_{P'_2},\ldots,h_{P'_2},h_{P_2},\ldots,h_{P_2},h_{\alpha^2_{\epsilon-r+1}},\ldots,h_{\alpha^2_{d}}\}\subseteq \{h_{\alpha_{1}^2},\ldots,h_{\alpha_{d}^2}\},$$
  	where  the multidegrees  of $h_{P'_2}$ and $h_{P_2}$ in the middle set are $\epsilon-i_0$ and $r$ respectively. Hence, $\tau_2 \subseteq \{h_{\alpha_{1}^2},\ldots,h_{\alpha_{d}^2}\}$.
  \end{proof}

  \begin{Remark}
  	(1)\label{replacingvertex} Let $p\in\{2,3\}$ and $h = \sum_{i=1}^{p+2}h_{\alpha_{i}}\in H$. Assume that there is an  $\alpha^2_0\in \MI(P_2')\setminus \{{\alpha_{1}^2},\ldots,{\alpha_{p+2}^2}\}$. Let $\widetilde{h} = \sum_{i=1}^{p+2}h_{\beta_{i}}$ be  an element of $H$ such that
  	\[ 
  	\beta_{i} =
  	\begin{cases} 
  	\alpha_{i} &\text{if} \ x\notin\alpha_{i} \text{,}  \\
  	\alpha^1_{i}\cup\alpha^2_0 & \text{if}\ x\in  \alpha_{i}.
  	\end{cases}
  	\] 
  	For  $\tau=\{h_{\gamma_{1}},\ldots, h_{\gamma_{m}}\}$,  define $\tau' :=\{h_{\nu_{1}},\ldots, h_{\nu_{m}}\}$, where 
  	\[ 
  	\nu_{j} =
  	\begin{cases} 
  	\gamma_{j} &\text{if} \ x\notin\gamma_{j} \text{,}  \\
  	\gamma^1_{j}\cup\alpha^2_0 & \text{if}\ x\in  \gamma_{i}.
  	\end{cases}
  	\] 
  	Then $\tau$ is a simplex of $\Delta_h$ if and only if $\tau'$ is a simplex of $\Delta_{\widetilde{h}}$. Therefore, $\Delta_h\cong \Delta_{\widetilde{h}}$.
  	
  	(2) Let $p\in\{2,3\}$ and $h = \sum_{i=1}^{p+2}h_{\alpha_{i}}\in H$. Let $\alpha^2,\widetilde{\alpha}^2 \in \{{\alpha_{1}^2},\ldots,{\alpha_{p+2}^2}\}$ with  $\alpha^2\neq\widetilde{\alpha}^2$. Let $\widetilde{h} = \sum_{i=1}^{p+2}h_{\beta_{i}}$ be  an element of $H$ such that 
  	\[
  	\beta_{i} =
  	\begin{cases} 
  	\alpha_{i} &\text{if} \quad \alpha^2_i\neq\alpha^2,\widetilde{\alpha}^2\text{,}  \\
  	\alpha^1_i\cup\widetilde{\alpha}^2 & \text{if}\quad  \alpha_i=\alpha^1_{i}\cup{\alpha}^2 \ \text{,}  \\
  	\alpha^1_i\cup{\alpha}^2 & \text{if}\quad \alpha_i=\alpha^1_{i}\cup\widetilde{\alpha}^2.
  	\end{cases}
  	\]
  	For  $\tau=\{h_{\gamma_{1}},\ldots, h_{\gamma_{m}}\}$,  define $\tau' :=\{h_{\nu_{1}},\ldots, h_{\nu_{m}}\}$, where 
  	\[
  	\nu_j =
  	\begin{cases} 
  	\gamma_j &\text{if} \quad \gamma^2_j\neq\alpha^2,\widetilde{\alpha}^2\text{,}  \\
  	\gamma^1_j\cup\widetilde{\alpha}^2 & \text{if}\quad  \gamma_j=\gamma^1_{j}\cup{\alpha}^2 \ \text{,}  \\
  	\gamma^1_j\cup{\alpha}^2 & \text{if}\quad \gamma_j=\gamma^1_{j}\cup\widetilde{\alpha}^2.
  	\end{cases}
  	\]
  	Observe that $\tau$ is a simplex of $\Delta_h$ if and only if $\tau'$ is a simplex of $\Delta_{\widetilde{h}}$. Therefore, $\Delta_h\cong \Delta_{\widetilde{h}}$.
  \end{Remark}
  
  \begin{Proposition}\label{homologoustozero}
  	Let  $h = (1,3,h_2, h_1)= \sum_{i=1}^{4}h_{\alpha_{i}}\in H$. Assume that $\MI(P_2')\subseteq\{{\alpha_{1}^2},\ldots,{\alpha_{4}^2}\}$ and there are exactly two $i$'s with $\alpha^2_i= P_2$. Let $\gamma$ be an 1-cycle in $\Delta_h$. If $\gamma$ is homologous to 0 in $X_h$, then it is also homologous to 0 in $\Delta_h$.
  \end{Proposition}
  \begin{proof}
  	Let $\eta=\sum_{}^{}c_{\sigma}\sigma$, where $c_{\sigma}\in \ZZ$, be a 2-chain in $X_h$ such that $\partial\eta=\gamma$. We construct an $\eta'$ in $\Delta_h$ such that $\partial\eta'=\gamma$. Let $\{h_{\nu_{1}},h_{\nu_{2}},h_{\nu_{3}}\}$ be a simplex in $\eta$. By Lemma~\ref{counting}, it is not a simplex of $\Delta_h$ if and only if  $\nu^2_j= P_2$ for $j=1,2,3.$  Let $\sigma = \{h_{\nu_{1}},h_{\nu_{2}},h_{\nu_{3}}\}$ be a simplex in $\eta$ such that it is not a simplex of $\Delta_h$. Note that $(0,4,h_2, h_1)-\sum_{i=1}^{3}h_{\nu_{i}}\in H$, call it $h_{\nu_{4}}$. Observe that $\{h_{\nu_{1}},\ldots,h_{\nu_{4}}\}$ is a face of $X_h$.  Define $$\sigma' := \sum_{j=1}^{3}(-1)^{j-1} \big\{ h_{\nu_{4}},h_{\nu_{1}},\ldots, \widehat{h_{\nu_{j}}},\ldots,h_{\nu_{3}}\big\}.$$ By Lemma~\ref{counting}, $\nu^2_4\neq P_2$. Therefore, $\sigma'$ is a 2-chain in $\Delta_h$. Observe that $\partial\sigma=\partial\sigma'$. Take $\eta'=\sum_{\sigma\in \eta}^{}c_{\sigma}\sigma'$, where $\sigma'=\sigma$ if $\sigma\in \Delta_h$ otherwise $\sigma'$ is as defined above for $\sigma$.  Then $\eta'$ is a 2-chain in $\Delta_h$ and $\partial\eta'=\gamma$. This completes the proof.
  \end{proof}

  \begin{Remark}\label{subsimplicial}
  	Let $Q$ be a poset such that it is a disjoint union of a poset $P_1$ and a chain $P_2'=\{a_1,a_2,x\}$, with $a_1\lessdot a_2\lessdot x$. Assume that $R[\MI(P_1)]$ satisfies property $N_2$. Let $H$ be the semigroup associated to $R[\MI(Q)]$. Let $A=\{\beta^1_1,\ldots,\beta^1_4\}, B=\{\beta^1_1,\beta^1_2,\delta^1_1,\delta^1_2\}$ where $\beta^1_j, \delta^1_i\in \MI(P_1)$, be two multisets with $\{\delta^1_1,\delta^1_2\}\cap\{\beta^1_3,\beta^1_4\}=\emptyset$, $\beta^1_1\neq\beta^1_2$ and $\sum_{\beta\in A}^{}h_{\beta}=\sum_{\beta\in B}^{}h_{\beta}\in  H_{1}$. Let 
  	$$\mathcal{S}=\big\{\{\nu_1,\ldots,\nu_4\}\subseteq\MI(Q): \{\nu^2_1,\ldots,\nu^2_4\}=\MI(P'_2)\ \text{and}\  \{\nu^1_1,\ldots,\nu^1_4\} \in\{A,B\} \big\}.$$
  	Let $\Delta'$ be the simplicial  complex whose facets are $\{h_{\nu_{1}},\ldots,h_{\nu_{4}}\}$, where $\{\nu_1,\ldots,\nu_4\}\in \mathcal{S}$. The choices of $A$ and $B$, up to isomorphism, are following:\par 
  	$(a)$ $\{\beta^1_1,\beta^1_2\}=\{\beta^1_3,\beta^1_4\}$, $\delta^1_3\neq\delta^1_4$,\par
  	$(b)$ $\{\beta^1_1,\beta^1_2\}=\{\beta^1_3,\beta^1_4\}$, $\delta^1_3=\delta^1_4$,\par
  	$(c)$ $\beta^1_1=\beta^1_3$, $\beta^1_2\neq\beta^1_4$, $\delta^1_3\neq\delta^1_4$ and $\beta^1_2\notin\{\delta^1_3,\delta^1_4\},$\par
  	$(d)$ $\beta^1_1=\beta^1_3$, $\beta^1_2\neq\beta^1_4$,$\delta^1_3$ and $\delta^1_3=\delta^1_4$,\par
  	$(f)$ $\beta^1_1=\beta^1_3$, $\beta^1_2\neq\beta^1_4$, $\delta^1_3\neq\delta^1_4$ and $\beta^1_2=\delta^1_3$, \par
  	$(g)$ $\beta^1_1=\beta^1_3$, $\beta^1_2\neq\beta^1_4$ and $\beta^1_2=\delta^1_3=\delta^1_4$, \par	
  	$(h)$ $\beta^1_1=\beta^1_3=\beta^1_4$,  $\delta^1_3\neq\delta^1_4$ and $\beta^1_2\notin\{\delta^1_3,\delta^1_4\}$, \par	
  	$(i)$ $\beta^1_1=\beta^1_3=\beta^1_4$, $\beta^1_2\neq\delta^1_3$ and $\delta^1_3=\delta^1_4$,\par	
  	$(j)$ $\beta^1_1=\beta^1_3=\beta^1_4$, $\delta^1_3\neq\delta^1_4$ and $\beta^1_2=\delta^1_3$, \par
  	$(k)$ $\beta^1_1=\beta^1_3=\beta^1_4$, $\beta^1_2=\delta^1_3=\delta^1_4$, \par	
  	$(l)$ $\{\beta^1_1,\beta^1_2\}\cap\{\beta^1_3,\beta^1_4\}=\emptyset$, $\beta^1_3=\beta^1_4$, $\delta^1_3=\delta^1_4$, \par	
  	$(m)$ each element of $A$ and $B$ appears with multiplicity 1 and $A\cap B = \{\beta^1_1,\beta^1_2\}$.\\
  	One can use a computer to check that $\widetilde{H}_{1}(\Delta')= 0$.
  \end{Remark}
  
  \begin{Proposition}\label{n=2homologous}
  	Let  $h = (1,3,h_2, h_1)= \sum_{i=1}^{4}h_{\alpha_{i}}\in H$. Assume that $\{{\alpha_{1}^2},\ldots,{\alpha_{4}^2}\}=\MI(P_2')$ (as a multiset). Let $\gamma$ be an 1-cycle in $\Delta_h$. If $\gamma$ is homologous to 0 in $X_h$, then it is also homologous to 0 in $\Delta_h$.
  \end{Proposition}
  \begin{proof}
  	Let $\eta=\sum_{}^{}c_{\sigma}\sigma$, where $c_{\sigma}\in \ZZ$ be a 2-chain in $X_h$ such that $\partial\eta=\gamma$. Let $\{h_{\nu_{1}},h_{\nu_{2}},h_{\nu_{3}}\}$ be a simplex in $\eta$. By Lemma~\ref{counting}, it is not a simplex of $\Delta_h$ if and only if there exist exactly two $j's$ with $\nu^2_j=  P_2$. We prove the proposition by induction on $k^\eta := \sum|c_{\sigma}|,$ where $\sigma \ \text{is a simplex of} \ \eta \ \text{but it is not a simplex of}\ \Delta_h $.
  	
  	Let $\sigma_1 = \{h_{\nu_{1}},h_{\nu_{2}},h_{\nu_{3}}\}$ be a simplex in $\eta$ such that it is not a simplex of $\Delta_h$ and $\nu_{j_{1}},\nu_{j_{2}} =P_2$. Let  $a$ be the sign of the coefficient of $\sigma_1$ in $\eta$. Observe that $\{h_{\nu_{j_{1}}},h_{\nu_{j_{2}}}\} $ is a simplex in $\partial\sigma_1$ and it is not in $\Delta_h$. Since $\partial\eta$ is in $\Delta_h$, there is another $\{h_{\nu_{j_{1}}},h_{\nu_{j_{2}}}\}$ in $\partial\eta$ with the  opposite sign. So there is a simplex $\sigma_2$ of $\eta$ but not a simplex of $\Delta_h$ such that $\sigma_1\neq\sigma_2$ and  $\partial(a\sigma_1+b\sigma_2)$ is an 1-cycle in $\Delta_h,$
  	where $b$ is the sign of the coefficient of $\sigma_2$ in $\eta$.
  	
  	Let $\sigma_1,\sigma_2$ be as above. We will define an $\eta_1$ such that $k^{\eta_1}< k^{\eta}$. Suppose $\sigma_1$ and $\sigma_2$ are the faces of the same facet, say $F$.  Let $\sigma_3$ and $\sigma_4$ be other two faces of $F$. Then $\sigma_3, \sigma_4\in\Delta_h$, by Lemma~\ref{counting} and there exist $c,d\in \{1,-1\}$ such that $\partial(c\sigma_3+d\sigma_4)=\partial(a\sigma_1+b\sigma_2)$. Define  
  	$$\eta_1 := \eta-(a\sigma_1+b\sigma_2)-(c\sigma_3+d\sigma_4).$$
  	Observe that $\partial\eta_1=\partial\eta$. \par
  	
  	On the other hand, suppose $\sigma_1$ and $\sigma_2$ are not the faces of the same facet. For $i=1,2$, let $F_i$ be the facet of $X_h$ such that $\sigma_i$ is a face of $F_i$. Write $F_1 = \{h_{\beta_{1}},\ldots,h_{\beta_{4}}\}$ and $F_2 = \{h_{\beta_{1}},h_{\beta_{2}},h_{\delta_{1}},h_{\delta_{2}}\}$. Let $A=\{\beta^1_1,\ldots,\beta^1_4\}, B=\{\beta^1_{1},\beta^1_{2},\delta^1_1,\delta^1_2\}$. For $A$ and $B$, let $\Delta'$ be as defined in Remark~\ref{subsimplicial}. Observe that $\Delta'$ is a subsimplicial complex of $\Delta_h$ and $\partial(a\sigma_1+b\sigma_2)$ is an 1-cycle in $\Delta'$. Since $\widetilde{H}_{1}(\Delta')= 0$, there exists a 2-chain $\mu_{\sigma_1,\sigma_2}$ in $\Delta'$ such that $\partial\mu_{\sigma_1,\sigma_2}=\partial(a\sigma_1+b\sigma_2)$. Define
  $$\eta_1 := \eta-(a\sigma_1+b\sigma_2)-\mu_{\sigma_1,\sigma_2}.$$
  	Observe that $\partial\eta_{1}=\partial\eta$. Also, notice that in both cases, $k^{\eta_1}< k^{\eta}$. Hence the proof.
  \end{proof}

  \begin{proof}[Proof of Proposition~\ref{bettiunion}]
  	Let $H$ and $H'$ be the semigroup associated to $R[\MI(Q)]$ and $R[\MI(P)]$ respectively. To prove the theorem, by Proposition~\ref{bh} and Lemma~\ref{vanishingbetti}, it suffices to show that for all  $\epsilon \in\{0,\ldots,4\}$, if $h=(\epsilon,4-\epsilon,h_2, h_1) =\sum_{i=1}^{4}h_{\alpha_{i}}\in H$ then $\widetilde{H}_{1}(\Delta_h)= 0$. We prove the theorem in the following two cases: $\MI(P_2')\nsubseteq \{{\alpha_{1}^2},\ldots,{\alpha_{4}^2}\}$ and  $\MI(P_2')\subseteq \{{\alpha_{1}^2},\ldots,{\alpha_{4}^2}\}$. In particular, if $n\geq3$, then we always have $\MI(P_2')\nsubseteq \{{\alpha_{1}^2},\ldots,{\alpha_{4}^2}\}$. 
  	\begin{itemize}
  		\item[$(a)$] Assume that $\MI(P_2')\nsubseteq \{{\alpha_{1}^2},\ldots,{\alpha_{4}^2}\}$. If $P_2'\in \{{\alpha_{1}^2},\ldots,{\alpha_{4}^2}\}$, then by Remark~\ref{replacingvertex}(1), $\Delta_h\cong\Delta_{\widetilde{h}}$, where $\widetilde{h}$ is as defined in Remark~\ref{replacingvertex}(1). Observe that $\widetilde{h}=(0,4,\widetilde{h'})$, where $\widetilde{h'}\in H'$. We know that $\Delta_{\widetilde{h}}\cong\Delta_{\widetilde{h'}}$. By hypothesis, $\widetilde{H}_{1}(\Delta_{\widetilde{h'}})= 0$. Therefore, $\widetilde{H}_{1}(\Delta_{\widetilde{h}})= 0$. If $P_2'\notin \{{\alpha_{1}^2},\ldots,{\alpha_{4}^2}\}$, then $h=(0,4,h')$, where $h'\in H'$.  Thus, $\widetilde{H}_{1}(\Delta_{h})= 0$. 
  		
  		\item[$(b)$] Now we assume that $\MI(P_2')\subseteq \{{\alpha_{1}^2},\ldots,{\alpha_{4}^2}\}$. We prove this case in two subcases  $n=1$ and $n=2$. If $n=1$, then by Remark~\ref{replacingvertex}(2), it is enough to consider the subcase $h=(1,3,h_2, h_1)$ and there are exactly two $\alpha^2_i$'s with $\alpha^2_i= P_2$. Let $\gamma$ be an 1-cycle in $\Delta_h$. By Proposition~\ref{homologous}, there exists an 1-cycle $\gamma'$ of $\Delta_{(0,4,h_2,h_1)}$ such that $\gamma$ is homologous to $\gamma'$ in $X_h$. By hypothesis, $\gamma'$ is homologous to 0 in $\Delta_h$. Thus, by Proposition~\ref{homologoustozero}, $\gamma$ is homologous to 0 in $\Delta_h$. This concludes the proof for $n=1$. For $n=2$, if $\MI(P_2')\subseteq \{{\alpha_{1}^2},\ldots,{\alpha_{4}^2}\}$, then $\MI(P_2')= \{{\alpha_{1}^2},\ldots,{\alpha_{4}^2}\}$. By the similar argument of the case $n=1$ and Proposition~\ref{n=2homologous}, we are done in this subcase also. Hence the proof.	
  	\end{itemize} 
  \end{proof}

 \section{Minimal Koszul syzygies of Hibi rings}\label{koszulpairs}

 Let $S= K[x_1,\ldots,x_n]$ be a standard graded polynomial ring over a field $K$. Let $I = (f_1,\ldots,f_m)$ be a graded ideal in $S$. Let $\{e_1,\ldots,e_m\}$ be a basis of the free $S$-module $S^m$. Define a map  $\varphi: S^m \to S$ by $\varphi(e_i) = f_i.$
 Then, $\ker \varphi$ is the second syzygy module of $S/I$, denoted by $Syz_2(S/I)$. Let $f_i$ and $f_j$ be two distinct generators of $I$. Then the Koszul relation $f_ie_j - f_je_i$ belongs to $Syz_2(S/I)$. We say $f_i
 , f_j$ a {\em Koszul relation pair} if $f_ie_j - f_je_i$ is a minimal generator of $Syz_2(S/I)$.\par
 
 Let $L = \MI(P)$ be a distributive lattice. Let $R[L] = K[L]/I_L$ be the Hibi ring associated to $L$. Let $<$ be a total order on the variables of $K[L]$ with the property that $x_\alpha < x_\beta$ if $\alpha < \beta$ in $L$. Consider the reverse lexicographic order $<$ on $K[L]$ induced by this order of the variables.  Recall from Section~\ref{hibiringdefi}, we have 
 \[\ini_<(I_L)=(x_\alpha x_\beta: \alpha,\beta\in L \ \text{and}\ \alpha, \beta \text{ incomparable}).\]
 Let us define $D_2 := \{(\alpha,\beta) : \alpha,\beta \in L\ \text{and}\  \alpha, \beta \text{ incomparable}\}$. 
 
 \subsection{Syzygies of initial Hibi ideals}\label{regsection}
 
 Let $K$ be a field and $\Delta$ be a simplicial complex on a vertex set $V =\{v_1,\ldots,v_n\}$. Let $K[\Delta]$ be the Stanley-Reisner ring of the simplicial complex $\Delta$. We know that $K[\Delta]= S/I_{\Delta}$, where $S = K[x_1,\ldots,x_n]$ and $I_{\Delta} = \{x_{i_1} \cdots x_{i_r} :\  \{v_{i_1},\ldots,v_{i_r}\} \notin \Delta  \}$. Since $K[\Delta]$ is a $\ZZ{^n}$-graded $S$-module, it has a minimal $\ZZ{^n}$-graded free resolution. Let $W \subset V $; we set $\Delta_W = \{F \in \Delta: F \subset W \}$. It is clear that $\Delta_W$ is again a simplicial complex.
 
 Let $L =\MI(P)$ be a distributive lattice. We know that $L$ is a poset under the order $\alpha \leq \beta$ if there is a chain from $\alpha$ to $\beta$. Let $\Delta(L)$ be the order complex of $L$. We have $K[\Delta(L)]= K[L]/I_{\Delta(L)}$, where $K[L] = K[\{x_\alpha: \alpha \in L \}]$ and $I_{\Delta(L)} = (x_{\alpha_{{1}}}\cdots x_{\alpha_{{r}}} :\ \{\alpha_{{1}},\ldots,\alpha_{{r}}\}  \notin \Delta(L))$.\par

 \begin{Lemma}
 	$I_{\Delta(L)} =  \ini_<(I_L)$.
 \end{Lemma}
 \begin{proof}
 	If $\alpha,\beta \in Q$ such that $\alpha$ and $\beta$ are incomparable, then $\{\alpha,\beta \} \notin \Delta(L)$. Hence, $x_\alpha x_\beta \in I_{\Delta(L)}$.\par
 	
 	On the other hand, if $x_{\alpha_1}\cdots x_{\alpha_r} \in I_{\Delta(L)}$, then $\{\alpha_1,\ldots,\alpha_r\}$ is not a chain. So there exist $\alpha,\beta \in \{\alpha_1,\ldots,\alpha_r\}$ such that $\alpha$ and $\beta$ are incomparable. Hence, $x_{\alpha_1}\cdots x_{\alpha_r} \in (x_{\alpha}x_{\beta}) \subseteq \ini_<(I_L)$. This concludes the proof.
 \end{proof}
 
 \begin{figure}
 	\begin{subfigure}[t]{4cm}
 		\centering	
 		\begin{tikzpicture}[scale=1]
 		\draw[fill= white] ;
 		\filldraw[black] (0,0) circle (.5pt) node[anchor=north] {};
 		\filldraw[black] (0.5,0) circle (.5pt) node[anchor=north] {};
 		\filldraw[black] (1,0) circle (.5pt) node[anchor=north] {};
 		\filldraw[black] (1.5,0) circle (.5pt) node[anchor=north] {};    
 		\end{tikzpicture}
 		\caption{$\widetilde{H}_{1}(\Delta;K) = 0$}
 	\end{subfigure}
 	\quad
 	\begin{subfigure}[t]{4cm}
 		\centering	
 		\begin{tikzpicture}[scale=1]
 		\draw[fill= white] (0,0)--(0,.9);
 		\filldraw[black] (0,0) circle (.5pt) node[anchor=north] {};
 		\filldraw[black] (0.5,0) circle (.5pt) node[anchor=north] {};
 		\filldraw[black] (1,0) circle (.5pt) node[anchor=north] {};
 		\filldraw[black] (0,0.9) circle (.5pt) node[anchor=north] {}; 
 		\end{tikzpicture}
 		\caption{$\widetilde{H}_{1}(\Delta;K) = 0$}
 	\end{subfigure}
 	\quad
 	\begin{subfigure}[t]{4cm}
 		\centering	
 		\begin{tikzpicture}[scale=1]
 		\draw[fill= white] (0,0)--(0,.9)--(.75,0);
 		\filldraw[black] (0,0) circle (.5pt) node[anchor=north] {};
 		\filldraw[black] (0.75,0) circle (.5pt) node[anchor=north] {};
 		\filldraw[black] (1.2,0) circle (.5pt) node[anchor=north] {};
 		\filldraw[black] (0,0.9) circle (.5pt) node[anchor=north] {}; 
 		\end{tikzpicture}
 		\caption{$\widetilde{H}_{1}(\Delta;K) = 0$}
 	\end{subfigure}
 	\quad
 	\begin{subfigure}[t]{4cm}
 		\centering	
 		\begin{tikzpicture}[scale=1]
 		\draw[fill= white] (0,0)--(0,.9)--(.75,0) (0,.9)--(1.5,0);
 		\filldraw[black] (0,0) circle (.5pt) node[anchor=north] {};
 		\filldraw[black] (0.75,0) circle (.5pt) node[anchor=north] {};
 		\filldraw[black] (1.5,0) circle (.5pt) node[anchor=north] {};
 		\filldraw[black] (0,0.9) circle (.5pt) node[anchor=north] {}; 
 		\end{tikzpicture}
 		\caption{$\widetilde{H}_{1}(\Delta;K) = 0$}
 	\end{subfigure}
 	\quad	
 	\begin{subfigure}[t]{4cm}
 		\centering
 		\begin{tikzpicture}[scale=1]
 		\draw[fill= white] (0,1)--(0,0)(.8,1)--(.8,0);
 		\filldraw[black] (0,0) circle (.5pt) node[anchor=north] {};
 		\filldraw[black] (.8,0) circle (.5pt) node[anchor=north] {};
 		\filldraw[black] (0,1) circle (.5pt) node[anchor=north] {};
 		\filldraw[black] (.8,1) circle (.5pt) node[anchor=north] {}; 
 		\end{tikzpicture}
 		\caption{$\widetilde{H}_{1}(\Delta;K) = 0$}
 	\end{subfigure}
 	\quad 	
 	\begin{subfigure}[t]{4cm}
 		\centering
 		\begin{tikzpicture}[scale=1]
 		\draw[fill= white] (0,0)--(0,1)--(1,0)--(1,1);
 		\filldraw[black] (0,0) circle (.5pt) node[anchor=north] {};
 		\filldraw[black] (1,0) circle (.5pt) node[anchor=north] {};
 		\filldraw[black] (0,1) circle (.5pt) node[anchor=north] {};
 		\filldraw[black] (1,1) circle (.5pt) node[anchor=north] {}; 
 		\end{tikzpicture}
 		\caption{$\widetilde{H}_{1}(\Delta;K) = 0$}
 	\end{subfigure}
 	\quad
  	\begin{subfigure}[t]{4cm}
  		\centering	
  		\begin{tikzpicture}[scale=1]
  		\draw[fill= white] (0,0)--(0,1)--(1,0)--(1,1)--(0,0);
  		\filldraw[black] (0,0) circle (.5pt) node[anchor=north] {};
  		\filldraw[black] (1,0) circle (.5pt) node[anchor=north] {};
  		\filldraw[black] (0,1) circle (.5pt) node[anchor=north] {};
  		\filldraw[black] (1,1) circle (.5pt) node[anchor=north] {}; 
  		\end{tikzpicture}
  		\caption{$\widetilde{H}_{1}(\Delta;K) = K$}\label{subfig}
  	\end{subfigure}
  	\quad	
 	\begin{subfigure}[t]{4cm}
 		\centering
 		\begin{tikzpicture}[scale=1]
 		\draw[fill= white] (0,0)--(0,.75)--(0,1.5);
 		\filldraw[black] (0,0) circle (.5pt) node[anchor=north] {};
 		\filldraw[black] (.8,0) circle (.5pt) node[anchor=north] {};
 		\filldraw[black] (0,.75) circle (.5pt) node[anchor=north] {};
 		\filldraw[black] (0,1.5) circle (.5pt) node[anchor=north] {}; 
 		\end{tikzpicture}
 		\caption{$\widetilde{H}_{1}(\Delta;K) = 0$}
 	\end{subfigure}
 	\quad
 	\begin{subfigure}[t]{4cm}
 		\centering
 		\begin{tikzpicture}[scale=1]
 		\draw[fill= white] (0,0)--(0,.75)--(0,1.5)(0,.75)--(.75,0);
 		\filldraw[black] (0,0) circle (.5pt) node[anchor=north] {};
 		\filldraw[black] (.75,0) circle (.5pt) node[anchor=north] {};
 		\filldraw[black] (0,.75) circle (.5pt) node[anchor=north] {};
 		\filldraw[black] (0,1.5) circle (.5pt) node[anchor=north] {}; 
 		\end{tikzpicture}
 		\caption{$\widetilde{H}_{1}(\Delta;K) = 0$}
 	\end{subfigure}	
 	\quad
 	\begin{subfigure}[t]{4cm}
 		\centering
 		\begin{tikzpicture}[scale=1]
 		\draw[fill= white] (0,0)--(0,.75)--(0,1.5)(0,1.5)--(.75,0);
 		\filldraw[black] (0,0) circle (.5pt) node[anchor=north] {};
 		\filldraw[black] (.75,0) circle (.5pt) node[anchor=north] {};
 		\filldraw[black] (0,.75) circle (.5pt) node[anchor=north] {};
 		\filldraw[black] (0,1.5) circle (.5pt) node[anchor=north] {}; 
 		\end{tikzpicture}
 		\caption{$\widetilde{H}_{1}(\Delta;K) = 0$}
 	\end{subfigure}
 	\quad
 	\begin{subfigure}[t]{4cm}
 		\centering
 		\begin{tikzpicture}[scale=1]
 		\draw[fill= white] (0,0)--(0,.75)--(0,1.5)--(0,2.25);
 		\filldraw[black] (0,0) circle (.5pt) node[anchor=north] {};
 		\filldraw[black] (0,2.25) circle (.5pt) node[anchor=north] {};
 		\filldraw[black] (0,.75) circle (.5pt) node[anchor=north] {};
 		\filldraw[black] (0,1.5) circle (.5pt) node[anchor=north] {}; 
 		\end{tikzpicture}
 		\caption{$\widetilde{H}_{1}(\Delta;K) = 0$}
 	\end{subfigure}
 	\caption{}\label{fig38}
 \end{figure}
 
 \begin{Theorem}\label{syzinitial}
 Let $(\alpha_1,\beta_1),(\alpha_2,\beta_2) \in D_2$. Then $x_{\alpha_1}x_{\beta_1},$ $x_{\alpha_2}x_{\beta_2}$ is a Koszul relation pair of $K[L]/\ini_<(I_L)$ if and only if either $\alpha_2 \vee \beta_2 \leq \alpha_1 \wedge \beta_1$ or $\alpha_1 \vee \beta_1 \leq \alpha_2 \wedge \beta_2$.
 \end{Theorem}
 
 \begin{proof}	
 	Suppose $x_{\alpha_1}x_{\beta_1},$ $x_{\alpha_2}x_{\beta_2}$ is a Koszul relation pair of $K[\Delta(L)]$. Let $W = \{\alpha_1, \beta_1, \alpha_2, \beta_2\}$. Then, by \cite[Theorem\ 5.5.1]{[BH93]}, $\widetilde{H}_{1}(\Delta_W;K) \neq 0$. All possible subsets of $L$ with cardinality $4$ are listed in Figure~\ref{fig38}. For $W' \subset L$ with $\#W' =4$, one can check that $\widetilde{H}_{1}(\Delta_{W'};K) \neq 0$  only if $W'$ is as in Figure~\ref{subfig}. Hence, the forward part follows.
 	
 	For the converse part, suppose $(\alpha_1,\beta_1),(\alpha_2,\beta_2) \in D_2$. Without loss of generality, assume that $\alpha_1 \vee \beta_1 \leq \alpha_2 \wedge \beta_2$. Let $W = \{\alpha_1, \beta_1, \alpha_2, \beta_2\}$. It is easy to see that  	
 	\[ \widetilde{H}_j(\Delta_W;K) = \begin{cases} 
 	K & \text{for} \quad j = 1, \\
 	0 & \text{for} \quad j \neq 1. \\
 	\end{cases}
 	\]
 	So by \cite[Theorem\ 5.5.1]{[BH93]}, $x_{\alpha_1}x_{\beta_1},$ $x_{\alpha_2}x_{\beta_2}$ is a Koszul relation pair of $K[\Delta(L)]$. Hence the proof.
 \end{proof}

 \subsection{Syzygies of Hibi ideals} 
 
  Let $L = \MI(P)$ be a distributive lattice with $\#L=n$. Let $R[L] = K[L]/I_L$ be the Hibi ring associated to $L$. There exists a weight vector $w =(w_1,\ldots,w_n)$ with strictly positive integer coordinates such that $\ini_{<_w} (I_L) = \ini_<(I_L)$ \cite[Theorem\ 22.3]{[PEEVA10]}.  Consider the polynomial ring $\widetilde{K[L]} = K[L][t]$ and the integral weight vector $\widetilde{w} = (w_1,...,w_n,1)$. Let $f = \sum_{i}^{} c_il_i \in K[L]$, where $c_i \in K \setminus \{0\}$ and $l_i$ is a monomial in $K[L]$. Let $l$ be a monomial in $f$ such that $w(l) = max_i\{w(l_i)\}$. Define $\tilde{f} =  \sum_{i}^{} t^{w(l)-w(l_i)}c_il_i$.  If we grade $\widetilde{K[L]}$ by deg$(t) = 1$ and deg$(x_i) = w_i$ for all $i$, then $\tilde{f}$ is homogeneous. Note that the image of $\tilde{f}$ in $\widetilde{K[L]}/(t-1)$ is $f$, and its image in $\widetilde{K[L]}/(t)$ is $\ini_{<_w} (f)$. Define $\widetilde{I_L} := (\widetilde{f}\ | f \in I_L)$. We have 
 $$\widetilde{I_L} =(\widetilde{x_{\alpha}x_{\beta}- x_{\alpha\wedge\beta} x_{\alpha\vee\beta}}  :\alpha,\beta\in L\ \text{and} \ \alpha,\beta \text{ incomparable}).$$
 
 \begin{Observation}\label{koszulrelation}
Let $(\alpha_1,\beta_1),(\alpha_2,\beta_2) \in D_2$. From the proof of \cite[Theorem\ 22.9]{[PEEVA10]}, we have
\begin{enumerate}
	\item[$(a)$] If $\widetilde{x_{\alpha_1}x_{\beta_1}- x_{\alpha_1\wedge\beta_1}x_{\alpha_1\vee\beta_1}},$  $\widetilde{x_{\alpha_2}x_{\beta_2}- x_{\alpha_2\wedge\beta_2}x_{\alpha_2\vee\beta_2}}$ is a Koszul relation pair of $\widetilde{K[L]}/\widetilde{I_L}$, then $x_{\alpha_1}x_{\beta_1},$ $x_{\alpha_2}x_{\beta_2}$ is a Koszul relation pair of $K[L]/\ini_<(I_L)$.
	\item[$(b)$] If $\widetilde{x_{\alpha_1}x_{\beta_1}- x_{\alpha_1\wedge\beta_1}x_{\alpha_1\vee\beta_1}},$  $\widetilde{x_{\alpha_2}x_{\beta_2}- x_{\alpha_2\wedge\beta_2}x_{\alpha_2\vee\beta_2}}$ is not a Koszul relation pair of $\widetilde{K[L]}/\widetilde{I_L}$, then $x_{\alpha_1}x_{\beta_1}-x_{\alpha_1\wedge\beta_1}x_{\alpha_1\vee\beta_1},$ $x_{\alpha_2}x_{\beta_2}-x_{\alpha_2\wedge\beta_2}x_{\alpha_2\vee\beta_2}$ is not a Koszul relation pair of $R[L]$.
    \item[$(c)$] From $(a)$ and $(b)$, we obtain that if 
	$x_{\alpha_1}x_{\beta_1}-x_{\alpha_1\wedge\beta_1} x_{\alpha_1\vee\beta_1}$, $x_{\alpha_2}x_{\beta_2}-x_{\alpha_2\wedge\beta_2}x_{\alpha_2\vee\beta_2}$ is a Koszul relation pair of $R[L]$, then $x_{\alpha_1}x_{\beta_1},$ $x_{\alpha_2}x_{\beta_2}$ is a Koszul relation pair of $K[L]/\ini_<(I_L)$.
	
\end{enumerate}
 \end{Observation}

 \begin{Theorem}\label{koszulhibiideal}
 	Let $(\alpha_1,\beta_1),(\alpha_2,\beta_2)  \in D_2$. If 
 	$x_{\alpha_1}x_{\beta_1}-x_{\alpha_1\wedge\beta_1}x_{\alpha_1\vee\beta_1},$ $x_{\alpha_2}x_{\beta_2}-x_{\alpha_2\wedge\beta_2}x_{\alpha_2\vee\beta_2}$ is a Koszul relation pair of $R[L]$, then either $\alpha_2 \vee \beta_2 \leq \alpha_1 \wedge \beta_1$ or $\alpha_1 \vee \beta_1 \leq \alpha_2 \wedge \beta_2$.
 \end{Theorem}
 
 \begin{proof}
 	The proof follows from Observation~\ref{koszulrelation}$(c)$ and Theorem~\ref{syzinitial}.	
 \end{proof}


\end{document}